\documentclass[10pt]{article}
\usepackage{amstext}
\usepackage{amsfonts}
\usepackage{amssymb}
\usepackage{amsbsy}
\usepackage{latexsym}
\usepackage{xy}
\usepackage{hhline}
\xyoption{all}
\newcommand{\vtx}[1]{*+[o][F-]{\scriptscriptstyle #1}} 
\newcounter{num}[section] %

\newenvironment{theo}
{\refstepcounter{num}%
\bigskip\noindent{\bf Theorem~\arabic{section}.\arabic{num}. }\it}

\newenvironment{cor}
{\refstepcounter{num}%
\bigskip\noindent{\bf Corollary~\arabic{section}.\arabic{num}. }\it}

\newenvironment{lemma}
{\refstepcounter{num}%
\bigskip\noindent{\bf Lemma~\arabic{section}.\arabic{num}. }\it}

\newcommand{\example}
{\refstepcounter{num}%
\bigskip\noindent{\bf Example~\arabic{section}.\arabic{num}.}}

\newcommand{\remark}
{\refstepcounter{num}%
\bigskip\noindent{\bf Remark~\arabic{section}.\arabic{num}.}}

\newcommand{\conj}
{\refstepcounter{num}%
\bigskip\noindent{\bf Conjecture~\arabic{section}.\arabic{num}.}}

\newcommand{\definition}[1]
{\refstepcounter{num}%
\bigskip\noindent{\bf Definition~\arabic{section}.\arabic{num}}~({\it #1}).}

\newenvironment{proof}{\medskip\noindent{\it Proof. }}
{$\Box$ \bigskip}
\newenvironment{proof_without_dot}{\medskip\noindent{\it Proof }}
{$\Box$ \bigskip}

\newenvironment{eq}{\begin{equation}}{\end{equation}}

\newcommand{\Ref}[1]{(\ref{#1})}

\newcommand{\si}{\sigma}
\newcommand{\al}{\alpha}
\newcommand{\be}{\beta}
\newcommand{\ga}{\gamma}
\newcommand{\la}{\lambda}

\newcommand{\ov}[1]{\overline{#1}}
\newcommand{\un}[1]{{\underline{#1}} }

\newcommand{\tr}{\mathop{\rm tr}}

\newcommand{\mdeg}{\mathop{\rm mdeg}}
\newcommand{\diag}{\mathop{\rm diag}}
\newcommand{\Char}{\mathop{\rm char}}

\newcommand{\sign}{\mathop{\rm{sgn }}}

\newcommand{\algA}{\mathcal{A}}    
\newcommand{\algB}{\mathcal{B}}    

\newcommand{\M}{\mathcal{M}}  
\newcommand{\N}{\mathcal{N}} 
\newcommand{\T}{{\mathcal T}} 
\newcommand{\I}{{\mathcal I}} 
\newcommand{\dec}{{\mathop{\rm{dec }}}} 
\newcommand{\Lin}{\mathop{\rm Lin}} 
\newcommand{\Glue}{\mathop{\rm Lin^{-1}}}

\newcommand{\FF}{{\mathbb{F}}}   
\newcommand{\NN}{{\mathbb{N}}}
\newcommand{\ZZ}{{\mathbb{Z}}}   
\newcommand{\QQ}{{\mathbb{Q}}}

\newcommand{\DP}{{\rm DP} }
\newcommand{\bpf}{{\mathop{\rm{bpf }}}}

\newcommand{\Q}{\mathcal{Q}}    
\newcommand{\G}{\mathcal{G}}    

\newcommand{\Sp}{S\!p}

\newcommand{\prII}{p} 
\newcommand{\prI}{\tilde{p}} 
\newcommand{\hc}[2]{{\rm hc}_{#1}(#2)} 
\newcommand{\Hc}[2]{{\rm hc}^{#1}(#2)} 

\newcommand{\binom}[2]{\left(\!\!\begin{array}{c} #1\\ #2\\
\end{array} \!\!\right)} 
\newcommand{\silin}[1]{\si_{#1}^{\rm lin}} 

\newcommand{\loopR}[3]{%
\begin{picture}(20,0)(#1,#2)
\put(-2,1){\llap{$\scriptstyle #3$}} \put(11,3){\circle{20}} \put(20,6){\vector(1,-4){1}}
\end{picture}}
\newcommand{\loopL}[3]{%
\begin{picture}(20,0)(#1,#2)
\put(22,1){$\scriptstyle #3$} \put(9,3){\circle{20}} \put(0,6){\vector(-1,-4){1}}
\end{picture}}

\newcommand{\rectangle}[2]{
\begin{picture}(0,0)
\put(-#1,-#2){\line(1,0){#1}}\put(0,-#2){\line(1,0){#1}}
\put(-#1,#2){\line(1,0){#1}}\put(0,#2){\line(1,0){#1}}
\put(-#1,-#2){\line(0,1){#2}}\put(-#1,0){\line(0,1){#2}}
\put(#1,-#2){\line(0,1){#2}}\put(#1,0){\line(0,1){#2}}
\end{picture}}

\begin{document}
\title{Relations between $O(n)$-invariants of several matrices} 
 \author{
A.A. Lopatin \\ 
{\small\it $Fakult\ddot{a}t$ $f\ddot{u}r$ Mathematik, %
$Universit\ddot{a}t$ Bielefeld, %
Postfach 100131, 33501 Bielefeld, Germany}  \\
{\small\it Institute of Mathematics, SBRAS, Pevtsova street, 13, Omsk 644099, Russia} \\
{\small\it artem\underline{ }lopatin@yahoo.com} \\
}
\date{} 
\maketitle 

\begin{abstract} A linear group $G<GL(n)$ acts on $d$-tuples of $n\times n$ matrices by simultaneous conjugation. In [Adv.~Math.~19 (1976), 306--381] Procesi established generators and relations between them for $G$-invariants, where $G$ is $GL(n)$, $O(n)$, and $\Sp(n)$ and the characteristic of base field is zero. We continue generalization of the mentioned results to the case of positive characteristic originated by Donkin in [Invent.~Math.~110 (1992), 389--401]. We investigate relations between generators for  
$O(n)$-invariants.
\end{abstract}




\section{Introduction}\label{section_intro}

\subsection{Relations}

We work over an infinite field $\FF$ of arbitrary characteristic $\Char{\FF}$. All vector spaces, algebras, and modules are over $\FF$ unless otherwise stated. By algebra we always mean an associative algebra. 

We consider the ring of polynomials  %
$$R=R_{n}=\FF[x_{ij}(k)\,|\,1\leq i,j\leq n,\, 1\leq k\leq d],$$ 
in $n^2 d$ variables. It is convenient to organize these variables in so-called {\it generic matrices}
$$X_k=\left(\begin{array}{ccc}
x_{11}(k) & \cdots & x_{1n}(k)\\
\vdots & & \vdots \\
x_{n1}(k) & \cdots & x_{nn}(k)\\
\end{array}
\right).$$ %
Denote coefficients in the characteristic polynomial
of an $n\times n$ matrix $X$ by $\si_t(X)$, i.e., %
\begin{eq}\label{eq1_intro} 
\det(X+\la E)=\sum_{t=0}^{n} \la^{n-t}\si_t(X).
\end{eq}
So, $\si_0(X)=1$, $\si_1(X)=\tr(X)$ and $\si_n(X)=\det(X)$.

Let $G$ be a group from the list $GL(n)$, $O(n)$, $\Sp(n)$, where we assume that $\Char{\FF}\neq2$ in case $G=O(n)$. The algebra of {\it matrix $G$-invariants} $R^G$ is the subalgebra of $R$ generated by $\si_t(A)$, where $1\leq t\leq n$ and $A$ ranges over all monomials in 
\begin{enumerate}
\item[$\bullet$] $X_1,\ldots,X_d$, if $G=GL(n)$; 

\item[$\bullet$] $X_1,\ldots,X_d$, $X_1^T,\ldots,X_d^T$, if $G=O(n)$;

\item[$\bullet$] $X_1,\ldots,X_d$, $X_1^{\ast},\ldots,X_d^{\ast}$, if $G=\Sp(n)$. (Here $X^{\ast}$ stands for the symplectic transpose of $X$).
\end{enumerate}
Moreover, we can assume that $A$ is {\it primitive}, i.e., is not equal to the power of a shorter monomial.  If the characteristic of $\FF$ is zero, then it is enough to take traces  $\tr(A)$, where $A$ can be non-primitive, in order to obtain $R^G$ 

For a field of characteristic zero generators for matrix
invariants of $G\in\{GL(n),O(n),\Sp(n)\}$ were described by Sibirskii
in~\cite{Sibirskii_1968} and Procesi in~\cite{Procesi76}. Generators for $G=SO(n)$ were calculated by Aslaksen et al.~in~\cite{Aslaksen95}. Over a field of arbitrary characteristic generators for matrix $GL(n)$-invariants were found by Donkin in~\cite{Donkin92a}. In~\cite{Zubkov99} Zubkov described generators for matrix $O(n)$- and
$\Sp(n)$-invariants, where in the orthogonal case he assumed that $\Char{\FF}\neq2$. Generators for matrix $SO(n)$-invariants were calculated by Lopatin in~\cite{Lopatin_so_inv}.

In case of characteristic zero relations between generators for matrix $GL(n)$-,
$O(n)$- and $\Sp(n)$-invariants were established by Procesi in~\cite{Procesi76}.
Independently, relations for $R^{GL(n)}$ were found by Razmyslov 
in~\cite{Razmyslov74}. Over a field of arbitrary characteristic relations for matrix
$GL(n)$-invariants were established by Zubkov in~\cite{Zubkov96} (see also Theorem~\ref{theo_relations_GL}). 

In this paper we describe relations for matrix $O(n)$-invariants modulo free relations in case $\Char{\FF}\neq2$ (see Theorem~\ref{theo_relations}). As a corollary, we explicitly establish relations between indecomposable invariants (see Corollary~\ref{cor_free}), where an element is called {\it decomposable} if it belongs to $\FF$-span of products of homogeneous invariants of positive degree.  

%

The above mentioned systems of generators are infinite, but it can be showed that they contain finite systems of generators. Moreover, the goal of the constructive theory of invariants is
to find out a minimal (by inclusion) homogeneous system of
generators of $K[V]^G$ explicitly. It is an important
problem, which arose as early as the theory of invariants itself. To find out a minimal system of generators we need a description of relations only between indecomposable invariants. And Corollary~\ref{cor_free} gives us such description for $R^{O(n)}$. 
As an application, in the upcoming paper~\cite{Lopatin_O3} some estimations are given on the highest degree of elements of a minimal system of generators for $R^{O(3)}$ (see Theorem~\ref{theo_free} below).

For more information on finite generating systems for matrix invariants see overviews~\cite{Formanek_1987} and~\cite{Formanek_1991} by Formanek. For recent developments in characteristic zero case see~\cite{Drensky_survey_2007} and in positive characteristic case see~\cite{DKZ_2002}. For results concerning generators and relations between them for mixed representations of quivers see survey~\cite{LZ_survey}.

Note that it is possible to define $O(n)$ and $SO(n)$ in characteristic two  
case. But in this case even generators for invariants of several vectors are not known
(for the latest developments see~\cite{Domokos_Frenkel}).

\smallskip
We introduce the following notions.
\begin{enumerate}
\item[$\bullet$] Let $\M$ be the monoid (without unity) freely generated by {\it letters}  $1,\ldots,d,1^T,\ldots,d^T$. 

\item[$\bullet$] Let $\N\subset\M$  be the subset of primitive elements, where the notion of a primitive element is defined as above.

\item[$\bullet$] Let $\M_{\FF}$ be the vector space with the basis $\M$.
\end{enumerate}
Assume that $\al=\al_1\cdots \al_p$ and $\be$ are elements of $\M$, where $\al_1,\ldots,\al_p$ are letters. 
\begin{enumerate}
\item[$\bullet$] Introduce an involution on $\M$ as follows. Define $\be^{TT}=\be$ for a letter $\beta$ and $\al^T=\al_p^T\cdots \al_1^T\in\M$. We extend the introduced map to an involution ${}^T:\M_{\FF}\to\M_{\FF}$ by linearity.

\item[$\bullet$] We say that $\al$ and $\be$ are {\it equivalent} and write $\al\sim\be$
if there exists a cyclic permutation $\pi\in S_p$ such that
$\al_{\pi(1)}\cdots\al_{\pi(p)}=\be$ or $\al_{\pi(1)}\cdots\al_{\pi(p)}=\be^T$.

\item[$\bullet$] Let $\M_{\si}$ be a ring with unity of (commutative) polynomials over $\FF$ freely generated by ``symbolic'' elements $\si_t(\al)$, where $t>0$ and $\al\in\M_{\FF}$.   

\item[$\bullet$] Let $\N_{\si}$ be a ring with unity of (commutative) polynomials over $\FF$ freely generated by ``symbolic'' elements $\si_t(\al)$, where $t>0$ and $\al\in\N$ ranges over $\sim$-equivalence classes. We can also define $\N_{\si}$ as a factor of $\M_{\si}$ by some ideal (see Lemma~\ref{lemma1_definition}). In particular, we can consider $\si_t(\al)$ with $\al\in\M_{\FF}$ as an element of $\N_{\si}$. 
\end{enumerate}
We will use the following convention:
$$\tr(\al)=\si_1(\al)$$
for any $\al\in\M_{\FF}$.   For a letter $\beta\in\M$ define
$$X_{\beta}=
\left\{
\begin{array}{rl}
X_{i},&\text{if } \beta=i\\
X_{i}^T,&\text{if } \beta=i^T\\
\end{array}
\right..
$$
Given $\al=\al_1\cdots \al_p\in\M$, where $\al_i$ is a letter, we assume that $X_{\al}=X_{\al_1}\cdots X_{\al_p}$. 

Consider a surjective homomorphism 
$$\Psi_n:\N_{\si} \to R^{O(n)}$$ %
defined by $\si_t(\al) \to \si_t(X_{\al})$, if $t\leq n$, and $\si_t(\al) \to 0$ otherwise. Its kernel $K_{n}$ is the ideal of {\it relations} for $R^{O(n)}$. Elements of $\bigcap_{i>0} K_{i}$ are called {\it free} relations.  For $\al,\be,\ga\in\M_{\FF}$ and $t,r\in\NN$, an element $\si_{t,r}(\al,\be,\ga)\in\N_{\si}$ is defined in Section~\ref{section_definition} (see Definition~\ref{def1_definition}), where ${\NN}$ stands for the set of non-negative integers.

\begin{theo}\label{theo_relations} If $\Char{\FF}\neq2$, then the ideal of relations $K_{n}$ for $R^{O(n)}\simeq  \N_{\si}/K_{n}$ is generated by 
\begin{enumerate}
\item[(a)] $\si_{t,r}(\al,\be,\ga)$, where $t+2r>n$ and 
$\al,\be,\ga\in\M_{\FF}$;

\item[(b)] free relations.
\end{enumerate} 
\end{theo}
\smallskip
Since $\FF$ is infinite, in the formulation of Theorem~\ref{theo_relations} we can take the following elements
\begin{enumerate}
\item[$\rm (a')$] $\si_{\un{t},\un{r},\un{s}}(x_1,\ldots,x_u, y_1,\ldots,y_v, z_1,\ldots,z_w)$ for $t+2r>n$, where $\un{t}=(t_1,\ldots,t_u)\in\NN^u$, $\un{r}=(r_1,\ldots,r_v)\in\NN^v$, $\un{s}=(s_1,\ldots,s_w)\in\NN^w$, $t=t_1+\cdots+t_u$, $r=r_1+\cdots+r_v=s_1+\cdots+s_w$, and $x_1,\ldots,z_w\in\M$ 
\end{enumerate}
instead of (a), where $\si_{\un{t},\un{r},\un{s}}$ is given by Definition~\ref{def_part_si}.

In Lemma~\ref{lemma1_free} we show that the ideal of free relations lies in the ideal of $\N_{\si}$ generated by $\si_t(\al)^{\Char{\FF}}$, where $\al$ ranges over $\N$ and $t>0$.  The following conjecture is valid in characteristic zero case.

\conj\label{conj1} In case $G=O(n)$ there are no non-zero free relations.
\bigskip

For $f\in R$ denote by $\mdeg{f}$ its {\it multidegree}, i.e.,
$\mdeg{f}=(t_1,\ldots,t_d)\in\NN^d$, where $t_k$ is the total degree of the polynomial $f$ in $x_{ij}(k)$, $1\leq i,j\leq n$. The algebra $R^{O(n)}$ is homogeneous with respect to $\NN^d$-grading. Denote the degree of $\al\in\M$ by $\deg{\al}$, the degree of $\al$ in a letter $\be$ by $\deg_{\be}{\al}$, and the multidegree of $\al$ by 
$$\mdeg(\al)=%
(\deg_{1}{\al}+\deg_{1^T}{\al},\ldots,\deg_{d}{\al}+\deg_{d^T}{\al}).$$
For $t>0$ we assume that $\deg{\si_t(\al)}=t\deg{\al}$ and $\mdeg{\si_t(\al)}=t\mdeg{\al}$. 
Therefore $\N$, $\M$, and $\N_{\si}$ have $\NN$-gradings by degrees and $\NN^d$-gradings by multidegrees.  Note that the given above generating system for $R^{O(n)}$ as well as relations $\rm (a')$ are $\NN^d$-homogeneous.

\subsection{The known results on relations}\label{section_PR}

In this section we explicitly formulate the results on relations that have been mentioned in Section~\ref{section_intro}. Let us remark that modulo free relations, Theorem~\ref{theo_relations} generalizes Theorem~\ref{theo_P} in the same manner as Theorem~\ref{theo_relations_GL} generalizes Theorem~\ref{theo_P_R}.

\begin{theo}(Procesi~\cite{Procesi76})\label{theo_P} If $\Char{\FF}=0$, then the ideal of relations $K_{n}$ for $R^{O(n)}\simeq  \N_{\si}/K_{n}$ is generated by $\si_{t,r}(\al,\be,\ga)$, where $t+2r=n+1$ and $\al,\be,\ga\in\M_{\FF}$.
\end{theo}
\bigskip

We denote by $\M'$  the monoid (without unity) freely generated by letters  $1,\ldots,d$ and denote by $\N'\subset\M'$ the subset of primitive elements. We define $\M'_{\FF}$ ($\M'_{\si}$, respectively) similarly to $\M_{\FF}$ ($\M_{\si}$, respectively). Let $\N'_{\si}$ be a ring of (commutative) polynomials over $\FF$ freely generated by ``symbolic'' elements $\si_t(\al)$, where $t>0$ and $\al\in\N'$ ranges over all primitive {\it cycles} (i.e., equivalence classes with respect to cyclic permutations). An analogue of Lemma~\ref{lemma1_definition} is also valid for $\N'_{\si}$. Hence we can consider $\si_{t}(\al)$ with $\al\in\M'_{\FF}$ as an element of $\N'_{\si}$.

\begin{theo}(Razmyslov~\cite{Razmyslov74}, Procesi~\cite{Procesi76})\label{theo_P_R} If $\Char{\FF}=0$, then the ideal of relations $K'_{n}$ for $R^{GL(n)}\simeq  \N'_{\si}/K'_{n}$ is generated by $\si_{n+1}(\al)$, where $\al\in\M'_{\FF}$.
\end{theo}

\begin{theo}(Zubkov,~\cite{Zubkov96})\label{theo_relations_GL} The ideal of relations $K'_{n}$ for $R^{GL(n)}\simeq  \N'_{\si}/K'_{n}$ is generated by $\si_{t}(\al)$, where $t>n$ and $\al\in\M'_{\FF}$.
\end{theo}

\remark\label{rem21} In case $G=GL(n)$ there are no non-zero free relations (see~\cite{Donkin93a}).
\bigskip

\subsection{Historical remarks on $\si_{t,r}$}

The complete linearization of $\si_{t,r}$ was introduced by
Procesi (see~\cite{Procesi76}, Section~8 of Part~I), where it was denoted by $F_{k,n+1}$. It is not difficult to see that 
\begin{eq}\label{eqPr}
\silin{t,r}(x_1,\ldots,x_t,y_1,\ldots,y_r,z_1,\ldots,z_r)= F_{r,t+2r}(y_1,\ldots,y_r,x_1,\ldots,x_t,z_1,\ldots,z_r)\text{ in }\N_{\si},
\end{eq}%
where $\silin{t,r}$ stands for the complete linearization of $\si_{t,r}$ (see Definition~\ref{def_part_si}) and $x_1,\ldots,z_r\in\M_{\FF}$. Then $\si_{t,r}$ was introduced by Zubkov in~\cite{ZubkovII}. Note that the definition from~\cite{ZubkovII} is different from our definition and their equivalence is established in Lemma~\ref{lemmaZ_appendix}. 

Another way to define $\si_{t,r}$ is via the determinant-pfaffian $\DP_{r,r}(X,Y,Z)$
that was introduced in~\cite{LZ1} as a ``mixture"{} of the determinant of $X$ and pfaffians
of $Y$ and $Z$ (see Example~\ref{ex1c_appendix}). Lemma~\ref{lemma5_si} and  Example~\ref{ex1c_appendix} imply that $\DP_{r,r}$ relates to $\si_{t,r}$ in the same way as the determinant relates
to $\si_t$, i.e., for $n\times n$ matrices $X,Y,Z$ we have
\begin{eq}\label{eq1_history}
\DP_{r,r}(X+\la E,Y,Z)=\sum_{t=0}^{t_0} \la^{t_0-t}\si_{t,r}(X,Y,Z),
\end{eq}%
where $n=t_0+2r$ and $t_0\geq0$. Since $\si_{t,0}(X,Y,Z)=\si_t(X)$ (see Remark~\ref{rem01}), Formula~\Ref{eq1_history} turns into~\Ref{eq1_intro} for $r=0$. 
Note that this
approach gives us $\si_{t,r}(X,Y,Z)$ as a polynomial in entries of matrices $X,Y,Z$. But
for our purposes we have to present $\si_{t,r}(X,Y,Z)$ in a different way, namely, as a
polynomial in $\si_t(\alpha)$, where $t$ ranges over positive integers and $\al$ ranges
over monomials in $X,Y,Z,X^T,Y^T,Z^T$.

\subsection{The structure of the paper}

The paper is organized as follows. In Section~\ref{section_notations} notations are given that are used throughout the paper. Key notion $\si_{t,r}$ is introduced in Section~\ref{section_definition}. 
In Section~\ref{section_si} we obtain some results on $\si_{t,r}$. In Lemma~\ref{lemma1_si} a formula for partial linearization of $\si_{t,r}$ is presented, which is similar to Formula~\Ref{eq1_definition} from the definition of $\si_{t,r}$. In characteristic zero case Lemma~\ref{lemma2a_si} allows us to work with complete linearization of $\si_{t,r}$ instead of $\si_{t,r}$ itself. 

In Section~\ref{section_proof} Theorem~\ref{theo_relations} is proven. The proof is based on results from paper~\cite{ZubkovII} by Zubkov, where an approach for computation of relations for $R^{O(n)}$ was given (see Sections~3 and~4 of~\cite{ZubkovII}). To complete Zubkov's approach we used the decomposition formula from~\cite{Lopatin_bplp}, which can be considered as a generalization of Amitsur's formula for $\si_{t,r}$.

In Section~\ref{section_free} we establish some restrictions on free relations (see Lemma~\ref{lemma1_free}). In particular, we show that there is no non-zero linear free relations between indecomposable invariants (see Corollary~\ref{cor_free}). 

In Section~\ref{section_appendix} the decomposition formula is formulated (see~\Ref{eq_decomp_appendix}). Then it is shown that our definition of $\si_{t,r}$ coincides with the original definition from~\cite{ZubkovII} (see Lemma~\ref{lemmaZ_appendix}). Since Section~\ref{section_appendix} contains numerous notions that have limited usage in the paper, we put it at the end of the paper.

\section{Notations}\label{section_notations}

In what follows, $\QQ$ stands for the quotient field of the ring of integers $\ZZ$. For a vector $\un{t}=(t_1,\ldots,t_u)\in\NN^u$ we write $\un{t}!=t_1!\cdots t_u!$. The length of $p$-tuple $\un{c}=(c_1,\ldots,c_p)$ we denote  by $\#\un{c}=p$. For short, we write $1^p$ for $(1,\ldots,1)\in\NN^p$. We assume that $E_m$ stands for the identity $m\times m$ matrix.

Consider $n,k\geq0$. If $k\leq n$, then we write $\binom{n}{k}$ for the binomial coefficient; otherwise, we define $\binom{n}{k}=0$. 
We also define $\binom{n}{-1}=0$.  Let us remark that $\binom{n}{0}=1$.

By substitution $1\to \al_1,\ldots,d\to \al_d$ in $\al\in\M$, where $\al_1,\ldots,\al_d\in\M$, we mean the substitution $1\to \al_1,\ldots,d\to \al_d$, $1^T\to \al_1^T,\ldots,d^T\to \al_d^T$ and denote it by $\al|_{1\to \al_1,\ldots,d\to \al_d}$. Similar convention we also use for substitutions of elements of $\N_{\si}$ and so on.

Denote by $\algB=\FF[x_1,\ldots,x_m]$ the polynomial ring in $x_1,\ldots,x_m$ over $\FF$, i.e., $\algB$ is a commutative $\FF$-algebra with unity generated by algebraically independent elements $x_1,\ldots,x_m$. 
%
Given an $\NN$-graded algebra $\algA=\sum_{i\geq0}\algA_i$ and $a=\sum_{i\geq 0}a_i\in\algA$, where almost all $a_i\in\algA_i$ are zero, we write $\Hc{s}{\algA}$ for $\sum_{0\leq i\leq s} A_i$ and $\hc{s}{a}=a_i$ for the homogeneous component of degree $s$.

A {\it quiver} $\Q=(\Q_0,\Q_1)$ is a finite oriented graph, where $\Q_0$ is the set of
vertices and $\Q_1$ is the set of arrows.  Multiple arrows and loops in $\Q$ are allowed.
For an arrow $\al$, denote by $\al'$ its head
and by $\al''$ its tail, i.e., %
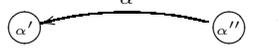
$$\vcenter{
\xymatrix@C=1cm@R=1cm{ %
\vtx{\al'}\ar@/^/@{<-}[rr]^{\al} && \vtx{\al''}\\
}} \quad.
$$
We say that $\al=\al_1\cdots \al_s$ is a {\it path} in $\Q$ (where $\al_1,\ldots,\al_s\in
\Q_1$), if $\al_1''=\al_2',\ldots,\al_{s-1}''=\al_s'$, i.e.,
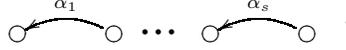
$$\vcenter{
\xymatrix@C=1cm@R=1cm{ %
\vtx{\;}\ar@/^/@{<-}[r]^{\al_1} & %
\vtx{\;} & %
\vtx{\;}\ar@/^/@{<-}[r]^{\al_s} & %
\vtx{\;}\\
}} \quad.
\begin{picture}(0,0)
\put(-73,-3){
\put(0,0){\circle*{2}} %
\put(-5,0){\circle*{2}} %
\put(5,0){\circle*{2}} %
} %
\end{picture}
$$
The head of the path $\al$ is $\al'=\al_1'$ and the tail is $\al''=\al_s''$.  A path
$\al$ is called {\it closed} if $\al'=\al''$. A closed path $\al$ is called {\it
incident} to a vertex $v\in\Q_0$ if $\al'=v$. Similarly, closed paths
$\be_1,\ldots,\be_s$ in $\Q$ are called {\it incident} to $v$ if $\be_1'=\cdots=\be_s'=v$.

\definition{of a mixed quiver} A quiver $\Q$ is called {\it mixed} if there are 
two maps ${}^T:\Q_0\to\Q_0$ and ${}^T:\Q_1\to\Q_1$ such that  
\begin{enumerate}
\item[$\bullet$] $v^{TT}=v$, $\be^{TT}=\be$;

\item[$\bullet$] $(\be^T)'=(\be'')^T$, $(\be^T)''=(\be')^T$
\end{enumerate} 
for all $v\in\Q_0$ and $\be\in\Q_1$.
\bigskip 

Assume that $\Q$ is a mixed quiver. Denote by $\M(\Q)$ the set of all closed paths in $\Q$ and denote by $\N(\Q)\subset \M(\Q)$ the subset of primitive paths.  Given a path $\al$ in $\Q$, we define the path $\al^{T}$ and introduce $\sim$-equivalence on $\M(\Q)$ in the same way as in Section~\ref{section_intro}. Moreover, we define $\M_{\FF}(\Q)$, $\M_{\si}(\Q)$, and $\N_{\si}(\Q)$ in the same way as $\M_{\FF}$, $\M_{\si}$, and $\N_{\si}$ have been defined in Section~\ref{section_intro}. The notions of degree and multidegree for elements of $\M(\Q)$ and $\N_{\si}(\Q)$ are introduced as above.   

\example\label{example_notations} Consider the quiver $\Q$:
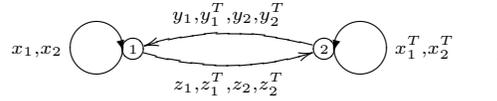
$$
\loopR{0}{0}{x_1,x_2} %
\xymatrix@C=1cm@R=1cm{ %
\vtx{1}\ar@1@/^/@{<-}[rr]^{y_1,y_1^T,y_2,y_2^T} &&\vtx{2}\ar@1@/^/@{<-}[ll]^{z_1,z_1^T,z_2,z_2^T}\\
}%
\loopL{0}{0}{x_1^T,x_2^T}\qquad\qquad,
$$
where there are four arrows from vertex $1$ to vertex $2$, there are four arrows in the
opposite direction, and each of its vertices has two loops. Here letters $x_i,x_i^T$, $y_i,y_i^T$, $z_i,z_i^T$ ($i=1,2$) stand for arrows of $\Q$. We define $1^T=2$ for vertex $1$, so $\Q$ is a mixed quiver. Then 
\begin{enumerate}
\item[$\bullet$] $(x_1 y_1^T z_1)^T=z_1^T y_1 x_1^T$, $y_1 z_1\sim y_1^T z_1^T$, $x_1 y_1 x_1^T z_1\sim x_1 y_1^T x_1^T z_1^T$;

\item[$\bullet$]  $\deg(y_1 z_1^T)=2$, $\deg_{y_1}(y_1 z_1)=1$, $\deg_{y_1^T}(y_1 z_1)=0$, and $\mdeg(x_1 x_1 x_2 y_1 x_1^T x_2^T z_2^T)=(3,2,\,1,0,\,0,1)$. 
\end{enumerate}

\section{The definition of $\si_{t,r}$}\label{section_definition}
In this section we assume that $\algA$ is a commutative unitary algebra over the field $\FF$ and all matrices are considered over $\algA$. 

Let us recall some formulas. In what follows $A,A_1,\ldots,A_p$ stand for $n\times n$ matrices and $1\leq t\leq n$. Amitsur's formula states~\cite{Amitsur_1980}:
\begin{eq}\label{eq_Amitsur}
\si_t(A_1+\cdots+A_p)=\sum (-1)^{t-(j_1+\cdots+j_q)} \si_{j_1}(\ga_1)\cdots\si_{j_q}(\ga_q),
\end{eq}
where the sum ranges over all pairwise different primitive cycles $\ga_1,\ldots,\ga_q$ in
letters $A_1,\ldots,A_p$ and positive integers $j_1,\ldots,j_q$ with
$\sum_{i=1}^{q}j_i\deg{\ga_i}=t$. Denote the right hand side of~\Ref{eq_Amitsur} by $F_{t,p}(A_1,\ldots,A_p)$. As an example,
\begin{eq}\label{eq_si2}
\si_2(A_1+A_2)=\si_2(A_1)+\si_2(A_2)+\tr(A_1)\tr(A_2)-\tr(A_1A_2).
\end{eq}
Note that for $a\in\algA$ we have 
\begin{eq}\label{eq20_definition}
\si_t(a A)=a^t\si_t(A).
\end{eq}

For $l\geq2$ we have the following well-known formula:
\begin{eq}\label{eq_D}
\si_t(A^l)=\sum\limits_{i_1,\ldots,i_{t l}\geq0}b^{(t,l)}_{i_1,\ldots,i_{t l}}
    \si_1(A)^{i_1}\cdots\si_{t l}(A)^{i_{t l}},%
\end{eq}
\noindent where we assume that $\si_i(A)=0$ for $i>n$. 
Denote the right hand side of~\Ref{eq_D} by $P_{t,l}(A)$. 
In~\Ref{eq_D} coefficients $b^{(t,l)}_{i_1,\ldots,i_{rl}} \in \ZZ$ do not depend on $A$
and $n$. If we take a diagonal matrix $A=\diag(a_1,\ldots,a_n)$, then $\si_t(A^l)$ is
a symmetric polynomial in $a_1,\ldots,a_n$ and $\si_i(A)$ is the
$i^{\rm th}$ elementary symmetric polynomial in $a_1,\ldots,a_n$, where $1\leq i\leq n$. Thus the coefficients
$b^{(t,l)}_{i_1,\ldots,i_{tl}}$ with $tl\leq n$ can easily be found. As an example,
\begin{eq}\label{eq_tr_a2}
\tr(A^2)=\tr(A)^2-2\si_2(A).
\end{eq}

\begin{lemma}\label{lemma1_definition} We have $\N_{\si}\simeq \M_{\si}/ L$ for the ideal $L$ generated by
\begin{enumerate}
\item[(a)] $\si_t(\al_1+\cdots+\al_p)-F_{t,p}(\al_1,\ldots,\al_p)$, 

\item[(b)] $\si_t(a\,\al)-a^t\si_t(\al)$,

\item[(c)] $\si_t(\al^l)-P_{t,l}(\al)$,

\item[(d)] $\si_t(\al\be)-\si_t(\be\al)$, 
 
\item[(e)] $\si_t(\al)-\si_t(\al^T)$,
\end{enumerate}
where $p>1$, $\al,\al_1,\ldots,\al_p\in \M_{\FF}$, $a\in\FF$, $t>0$, and $l>1$. 
\end{lemma}
\begin{proof} Consider a homomorphism $\rho:\M_{\si}\to \N_{\si}$ defined by $\si_t(\al)\to f$, where we apply equalities (a), (b), and (c) one after another to $\si_t(\al)$, $\al\in\M_{\FF}$, to obtain $f$. Obviously, $\rho$ is well defined. Since $\rho$ is surjective, to complete the proof it is enough to show that $\rho(f)=0$, where $f\in L$ ranges over elements (a)--(e). 

Assume that $f$ is not element (e) and $\rho(f)\neq0$. We define $\M'$, $\M'_{\si}$, and $\N'_{\si}$ in the same way as in Section~\ref{section_PR} with the only difference that now we start with the monoid $\M'$ freely generated by letters $1,\ldots,2d$. Note that $\M'\simeq\M$, where the isomorphism is given by $i\to i$ and $i+d\to i^T$ for all $1\leq i\leq d$. Hence $\M'_{\si}\simeq \M_{\si}$ and $f$ can be considered as an element of $\M'_{\si}$. Define $\rho':\M'_{\si}\to \N'_{\si}$ in the same way as $\rho$. Since $\rho(f)\neq 0$, we have $\rho'(f)\neq0$. But $\rho'(f)$ is a free relation for $R^{GL(n)}$; a contradiction to Remark~\ref{rem21}.

We assume that $f$ is element (e), i.e., %
$$f=\si_t(\sum_{i=1}^s a_i\al_i)-\si_t(\sum_{i=1}^s  a_i\al_i^T),$$ %
where $a_i\in\FF$ and $\al_i\in\M$. Denote the result of application of (a) and (b) to $f$ by $h\in\M_{\si}$. Let $a^{(i_1,\ldots,i_s)}$ stand for $a_1^{i_1}\cdots a_s^{i_s}$. Then
\begin{eq}\label{eq12}
h=\sum a^{j_1\mdeg{\ga_1}+\cdots+j_q\mdeg{\ga_q}}\,(-1)^{t-(j_1+\cdots+j_q)}\, %
(\si_{j_1}(\ga_1^{(1)})\cdots\si_{j_q}(\ga_q^{(1)})-
\si_{j_1}(\ga_1^{(2)})\cdots\si_{j_q}(\ga_q^{(2)})),
\end{eq}
where the sum ranges over $j_1,\ldots,j_q$ and cycles $\ga_1,\ldots,\ga_q$ in letters $x_1,\ldots,x_s$ as in~\Ref{eq_Amitsur}, 
$\ga_l^{(1)}$ is the result of substitution $x_1\to \al_1,\ldots, x_s\to\al_s$ in $\ga_l$,  and $\ga_l^{(2)}$ is the result of substitution  $x_1\to \al_1^T,\ldots, x_s\to\al_s^T$ in $\ga_l$ ($1\leq l\leq q$). Note that for any pair $(\un{j},\un{\ga})$ from~\Ref{eq12} there exists a pair $(\un{j},\un{\be})$ from~\Ref{eq12} such that $\ga_l^{(1)}\sim \be_l^{(2)}$ for all $1\leq l\leq q$. Applying (c) to $\si_{j_1}(\ga_1^{(1)})\cdots\si_{j_q}(\ga_q^{(1)})$ and $\si_{j_1}(\be_1^{(2)})\cdots\si_{j_q}(\be_q^{(2)})$, we obtain equal elements in $\N_{\si}$. Thus, $\rho(f)=0$.
\end{proof}


Let $t,r\in\NN$. In order to define $\si_{t,r}$, we consider the quiver $\Q$
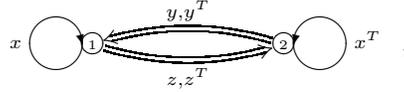
$$
\loopR{0}{0}{x} %
\xymatrix@C=1cm@R=1cm{ %
\vtx{1}\ar@2@/^/@{<-}[rr]^{y,y^T} &&\vtx{2}\ar@2@/^/@{<-}[ll]^{z,z^T}\\
}%
\loopL{0}{0}{x^T}\qquad,
$$
where there are two arrows from vertex $1$ to vertex $2$ and there are two arrows in the
opposite direction. We define $1^T=2$ for vertex $1$ to turn $\Q$ into a mixed quiver.

\definition{of $\si_{t,r}(x,y,z)$}\label{def1_definition} Denote by  ${\I}={\I}_{t,r}$ the set of pairs $(\un{j},\un{\al})$ such that 
\begin{enumerate}
 \item[$\bullet$] $\#\un{j}=\#\un{\al}=p$ for some $p$;

 \item[$\bullet$] $\al_1,\ldots,\al_p\in \N(\Q)$ belong to pairwise different $\sim$-equivalence classes and $j_1,\ldots,j_p\geq1$;

 \item[$\bullet$] $j_1\mdeg(\al_1)+\cdots+j_p\mdeg(\al_p)=(t,r,r)$.
\end{enumerate}
Then we define $\si_{t,r}(x,y,z)\in\N_{\si}(\Q)$ by
\begin{eq}\label{eq1_definition}
\si_{t,r}(x,y,z)=\sum_{(\un{j},\un{\al})\in{\I}} %
(-1)^{\xi} \;\si_{j_1}(\al_1)\cdots\si_{j_p}(\al_p), %
\end{eq}
where $p=\#\un{j}=\#\un{\al}$ and $\xi=\xi_{\un{j},\un{\al}}=t+\sum_{i=1}^p j_i(\deg_y{\al_i}+\deg_z{\al_i}+1)$. 
For $t=r=0$ we define $\si_{0,0}(x,y,z)=1$. Moreover, 
\begin{enumerate}
\item[$\bullet$] if $\al,\be,\ga\in\M_{\FF}$, then we define $\si_{t,r}(\al,\be,\ga)\in\N_{\si}$ as the result of substitution $x\to\al$, $y\to\be$, $z\to\ga$ in~\Ref{eq1_definition};

\item[$\bullet$] if $X,Y,Z$ are $n\times n$ matrices, then we define $\si_{t,r}(X,Y,Z)\in\algA$ as the result of substitution $x\to X$, $y\to Y$, $z\to Z$ in~\Ref{eq1_definition}, where we assume that $\si_j(A)=0$ for $t>n$ and any $n\times n$ matrix $A$.
\end{enumerate}

\example\label{ex_DP} {\bf 1.} If $t=0$ and $r=1$, then $\sim$-equivalence classes on $\N(\Q)$ are $y z,\, y z^T,\, \ldots$ %
Hence, $\si_{0,1}(x,y,z)=-\tr(yz)+\tr(yz^T)$.
\smallskip

{\bf 2.} If $t=r=1$, then $\sim$-equivalence classes on $\N(\Q)$ are 
$$x,\, y z,\, y z^T,\, x y z,\, x y z^T,\, x y^T z,\, x y^T z^T,\, \ldots$$
and we can see that $\si_{1,1}(x,y,z)=$
$$-\tr(x)\tr(yz)+\tr(x)\tr(yz^T)+ \\%
\tr(x y z)-\tr(x y z^T)-\tr(x y^T z)+\tr(x y^T z^T).$$ %

\remark\label{rem01} $\si_{t,0}(x,y,z)=\si_t(x)$.
\bigskip

\section{Properties of $\si_{t,r}$}\label{section_si}

In this section all matrices have entries in a commutative unitary algebra $\algA$ over the field $\FF$.

\definition{of a partial linearization of $\si_{t,r}$}\label{def_si}\label{def_part_si} We assume that $\un{t}=(t_1,\ldots,t_u)\in\NN^u$, $\un{r}=(r_1,\ldots,r_v)\in\NN^v$, $\un{s}=(s_1,\ldots,s_w)\in\NN^w$ satisfy
$$s_1+\cdots+s_v=r_1+\cdots+r_w,$$
and $x_1,\ldots,x_u$, $y_1,\ldots,y_v$, $z_1,\ldots,z_w$ belong to $\M_{\FF}$. Consider   
$\si_{t,r}(a_1 x_1+\cdots+a_u x_u, b_1 y_1+\cdots+b_v y_v, c_1 z_1+\cdots+c_w z_w)\in\N_{\si}$
as a polynomial in $a_1,\ldots,a_u$, $b_1,\ldots,b_v$, $c_1,\ldots,c_w$, where $a_1,\ldots, c_w\in\FF$. We denote the coefficient of $a_1^{t_1}\cdots a_u^{t_u}
b_{1}^{r_1}\cdots b_v^{r_v} c_{1}^{s_1}\cdots c_w^{s_w}$ in this polynomial by 
\begin{eq}\label{eq00_si}
\si_{\un{t},\un{r},\un{s}}(x_1,\ldots,x_u,y_1,\ldots,y_v,z_1,\ldots,z_w)\in\N_{\si}. %
\end{eq}
In other words, this coefficient is equal to the homogeneous component of  $\si_{t,r}(x_1+\cdots+x_u,y_1+\cdots+y_v,z_1+\cdots+z_w)$ of multidegree $(\un{t},\un{r},\un{s})$. In multilinear case we have $v=w$ and write %
$$\silin{u,v}(x_1,\ldots,x_u,y_1,\ldots,y_v,z_1,\ldots,z_w)$$ %
for $\si_{1^u,1^v,1^w}(x_1,\ldots,x_u,y_1,\ldots,y_v,z_1,\ldots,z_w)$, where the definition of $1^u$ was given in Section~\ref{section_notations}. 

Given for $n\times n$ matrices $X_1,\ldots,Z_w$, we define $\si_{\un{t},\un{r},\un{s}}(X_1,\ldots,X_u,Y_1,\ldots,Y_v,Z_1,\ldots,Z_w)\in\algA$ as the result of substitution $x_1\to X_1,\ldots, z_w\to Z_w$ in element~\Ref{eq00_si} of $\N_{\si}$ and %
use the assumption that $\si_j(A)=0$ for $j>n$ and any $n\times n$ matrix $A$.  
\bigskip

In this section we will show that an analogue of Formula~\Ref{eq1_definition} from the definition of $\si_{t,r}$ is valid for $\si_{\un{t},\un{r},\un{s}}$ (see Lemma~\ref{lemma1_si} below). Generalize the construction from Section~\ref{section_definition} as follows. 
Using notation from Definition~\ref{def_part_si}, we define the mixed quiver $\Q_{\un{t},\un{r},\un{s}}$:
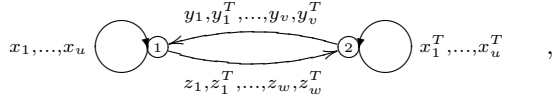
$$
\loopR{0}{0}{x_1,\ldots,x_u} %
\xymatrix@C=1cm@R=1cm{ %
\vtx{1}\ar@1@/^/@{<-}[rr]^{y_1,y_1^T,\ldots,y_v,y_v^T} &&\vtx{2}\ar@1@/^/@{<-}[ll]^{z_1,z_1^T,\ldots,z_w,z_w^T}\\
}%
\loopL{0}{0}{x_1^T,\ldots,x_u^T}\qquad\qquad\quad,
$$
where there are $2w$ arrows from vertex $1$ to vertex $2$, there are $2v$ arrows in the
opposite direction, and $1^T=2$. Here by abuse of notation $x_1,\ldots,z_w$ stand for the arrows of $\Q_{\un{t},\un{r},\un{s}}$ as well as for elements of $\M_{\FF}$ from the definition of $\si_{t,r}$. Hence we can assume that $\N_{\si}(\Q_{\un{t},\un{r},\un{s}})\subset \N_{\si}$.  

Given a path $\al$ of $\Q_{\un{t},\un{r},\un{s}}$, denote the total degree of $\al$ in $x_1,\ldots,x_u$ by $\deg_{x}{\al}$, and the total degree of $\al$ in $x_1^T,\ldots,x_u^T$ by $\deg_{x^T}{\al}$. Similar notation we also use for $y$ and $z$. 
As an example, $\deg_{x^T}(x_1 y_1 x_1^T x_2^T z_2^T)=2$. 

Assume that ${\I}={\I}_{\un{t},\un{r},\un{s}}$ is the set of pairs $(\un{j},\un{\al})$ such that 
\begin{enumerate}
 \item[$\bullet$] $\#\un{j}=\#\un{\al}=p$ for some $p$;

 \item[$\bullet$]  $\al_1,\ldots,\al_p\in \N(\Q_{\un{t},\un{r},\un{s}})$ are representatives of pairwise different $\sim$-equivalence classes and $j_1,\ldots,j_p\geq1$;

 \item[$\bullet$] $j_1\mdeg(\al_1)+\cdots+j_p\mdeg(\al_p)=(\un{t},\un{r},\un{s})$.
\end{enumerate}

\begin{lemma}\label{lemma1_si}  
We have  
$$\si_{\un{t},\un{r},\un{s}}(x_1,\ldots,x_u,y_1,\ldots,y_v,z_1,\ldots,z_w)= \sum_{(\un{j},\un{\al})\in {\I}} %
(-1)^{\eta} \;%
\si_{j_1}(\al_1)\cdots\si_{j_p}(\al_p),$$ %
where $p=\#\un{j}=\#\un{\al}$ and $\eta=\eta_{\un{j},\un{\al}}=t_1+\cdots+t_u+\sum_{i=1}^p j_i(\deg_y{\al_i}+\deg_z{\al_i}+1)$. 
\end{lemma}
\begin{proof} We set $t=t_1+\cdots+t_u$ and $r=r_1+\cdots+r_u=s_1+\cdots+s_w$. %
By definition, $\si_{\un{t},\un{r},\un{s}}(x_1,\ldots,x_u,y_1,\ldots,y_v,z_1,\ldots,z_w)$ is equal to the homogeneous component of multidegree $(\un{t},\un{r},\un{s})$ of 
$$f=\si_{t,r}(x_1+\cdots+x_u, y_1+\cdots+y_v, z_1+\cdots+z_w)$$
$$=\sum %
(-1)^{\xi} \;\si_{j_1}(\al_1)\cdots\si_{j_p}(\al_p)|_{x\to \sum_{i=1}^u x_i, y\to \sum_{j=1}^v y_j, z\to \sum_{k=1}^w  z_k},$$ %
where $\xi$, $\al_1,\ldots,\al_p$, and $j_1,\ldots,j_p$ are the same as in the definition of $\si_{t,s}$.  Let $\Q$ be the quiver from Section~\ref{section_definition}. For a closed path $\al$ in $\Q$ denote the result of substitution $x\to \sum x_i$, $y\to \sum y_j$, $z\to \sum z_k$ in $\al$ by %
$$\sum_{l\in L(\al)}\al^{(l)},$$
where $\al^{(l)}$ is a closed path in ${\Q}_{\un{t},\un{r},\un{s}}$ for all $l$. As an example, the result of substitution $x\to x_1+x_2$, $y\to y_1+y_2$, and $z\to z_1+z_2$ in $yz^T$ is equal to $y_1 z_1^T + y_1 z_2^T + y_2 z_1^T + y_2 z_2^T$.

Applying Amitsur's formula~\Ref{eq_Amitsur} to $f$, we obtain
\begin{eq}\label{eq_f}
f=\sum (-1)^{\xi}\cdot (-1)^{\xi_1} \si_{j_{11}}(\ga_{11})\cdots\si_{j_{1q_1}}(\ga_{1q_1}) %
\ast\cdots\ast  (-1)^{\xi_p} \si_{j_{p1}}(\ga_{p1})\cdots\si_{j_{pq_p}}(\ga_{pq_p}),
\end{eq}
where $\xi_i=j_i-(j_{i1}+\cdots + j_{iq_i})$. Here for any $1\leq i\leq p$ the sum ranges over  pairwise different primitive cycles $\ga_{i1},\ldots, \ga_{iq_i}$ in letters $\al_i^{(l)}$ (where $l\in L(\al_i)$) and $j_{i1},\ldots,j_{iq_i}\geq1$ such that $j_{i1}\deg(\ga_{i1})+\cdots+j_{iq}\deg(\ga_{iq_i})=j_i\deg(\al_i)$. 

We claim that taking the homogeneous component of $f$ of multidegree $(\un{t},\un{r},\un{s})$ we obtain the required formula for $\si_{\un{t},\un{r},\un{s}}$. We split the proof into several statements. 

\medskip 
\noindent{\it {\bf 1.} We have $(-1)^{\xi+\xi_1+\cdots+\xi_p}=
(-1)^{\eta}$, where $\eta=t+\sum_{i=1}^p \sum_{k=1}^{q_p} j_{ik}(\deg_y(\ga_{ik})+\deg_z(\ga_{ik})+1)$.}

\smallskip
\noindent {\it Proof.}
For $1\leq i\leq p$ we have $\deg_y{\al_i}=\deg_y{\al_i^{(l)}}$ for all $l$. Then 
$$\deg_y(\ga_{ik})=\deg_{y}{\al_i}\cdot \frac{\deg{\ga_{ik}}}{\deg{\al_i}}$$ 
for all $k$. Therefore
$$j_i\deg_y{\al_i}=\frac{\deg_y{\al_i}}{\deg{\al_i}} ( j_{i1}\deg(\ga_{i1})+\cdots+j_{iq}\deg(\ga_{iq_i}) )=$$
$$=j_{i1}\deg_y(\ga_{i1})+\cdots+j_{iq}\deg_y(\ga_{iq_i}).$$
The same formula is also valid for $z$. Now it is easy to see that the required formula is valid.

\medskip
Given a closed path $\ga$ in $\Q_{\un{t},\un{r},\un{s}}$, denote the result of substitution $x_1,\ldots,x_u\to x$, $y_1,\ldots,y_v\to y$, and $z_1,\ldots,z_w\to z$ in $\ga$ by $\psi(\ga)$. We can consider $\psi(\ga)$ as a closed path in $\Q$. Let us remark that if $\al^p=\be^q$ for $\al,\be\in \M(\Q)$, $p,q>0$, then there exists a $\ga\in \M(\Q)$ such that $\al=\ga^i$ and $\be=\ga^j$ for some $i,j$; in particular, if $\be$ is primitive, then $\al=\be^i$ for some $i$.

\medskip
\noindent{\it {\bf 2.} Let $\al$ be a primitive closed path in $\Q$ and $\ga$ be a primitive word in letters %
$\al^{(l)}$, where $l\in L(\al)$. Then $\ga$ is a primitive closed path in $\Q_{\un{t},\un{r},\un{s}}$. In particular, $\ga_{11},\ldots,\ga_{1q_1},\ldots,\ga_{p1},\ldots,\ga_{pq_p}$ are primitive closed paths. }

\smallskip
\noindent {\it Proof.} Assume that there exist $k>1$ and a primitive closed path $\be$ in $\Q_{\un{t},\un{r},\un{s}}$ such that $\ga=\be^k$. We have $\ga=\al^{(l_1)}\cdots\al^{(l_m)}$ for some $l_1,\ldots,l_m\in L(\al)$. Thus $\al^m=\psi(\be)^k$ and applying the above mentioned remark 
we obtain that there exists a $j>0$ such that $\psi(\be)=\al^{j}$. Since $\deg{\al^{(l)}}=\deg{\al}$ for all $l$, we have $\be=\al^{(l_1)}\cdots\al^{(l_{j})}$. Hence, %
$$\ga=(\al^{(l_1)}\cdots\al^{(l_{j})})^k=\al^{(l_1)}\cdots\al^{(l_m)}.$$  
Closed paths $\al^{(l)}$ (where $l\in L(\al)$) are pairwise different. So the last formula implies that $l_i=l_{i+j}$ for any $i\leq m-j$. In other words, $\ga$ is not a primitive word in letters $\al^{(l)}$, where $l\in L(\al)$; a contradiction.    
\smallskip

\medskip
\noindent{\it {\bf 3.} Closed paths $\ga_{11},\ldots,\ga_{1q_1},\ldots,\ga_{p1},\ldots,\ga_{pq_p}$ are pairwise different with respect to $\sim$-equivalence.}
    
\smallskip
\noindent {\it Proof.} Assume that $\ga_{ik}=\ga_{jm}$. Equalities $\psi(\ga_{ik})=\al_i\cdots \al_i$ and $\psi(\ga_{jm})=\al_{j}\cdots \al_{j}$ together with $\psi(\ga_{ik})=\psi(\ga_{jm})$ imply that $i=j$ (see the above mentioned remark). Since $\ga_{i1},\ldots,\ga_{iq_i}$ are pairwise different, we have $k=m$. 

\medskip
To complete the proof of the lemma we apply substitution $\psi$ in the same way as in the proof of Statement~3. 
\end{proof}

Denote by ${\I}^{\rm lin}_{u,v}$ the set of $\un{\al}$ such that $(1^p,\un{\al})\in{\I_{1^u,1^v,1^v}}$, where $p=\#\un{\al}$.
The following corollary is a trivial consequence of Lemma~\ref{lemma1_si}.

\begin{cor}\label{cor_si} We have
$$\silin{u,v}(x_1,\ldots,x_u,y_1,\ldots,y_v,z_1,\ldots,z_v)=\sum_{\un{\al}\in {\I}^{\rm lin}_{u,v}} %
(-1)^{\eta_{\un{\al}}} \tr(\un{\al}),$$
where $\eta_{\un{\al}}=u+\sum_{i=1}^p (\deg_y{\al_i}+\deg_z{\al_i}+1)$, $p=\#\un{\al}$, and $\tr(\un{\al})$ stands for $\tr(\al_1)\cdots \tr(\al_p)$.
\end{cor}

\begin{lemma}\label{lemma2a_si} 
Assume that $\FF=\QQ$. Then for $x,y,z\in\M_{\FF}$ we have 
$$\si_{t,r}(x,y,z)=\frac{1}{t!\,(r!)^2}\,\,
\silin{t,r}(\underbrace{x,\ldots,x}_{t}, \underbrace{y,\ldots,y}_{r}, 
    \underbrace{z,\ldots,z}_{r} ).$$
\end{lemma}
\begin{proof} In this proof we use notions formulated in Section~\ref{section_appendix} (see below).
Consider the tableau with substitution $(\T,(X,Y,Z))$ of dimension $(t+2r,t+2r)$ from Example~\ref{ex1_appendix}. By Lemma~\ref{lemma1_appendix}, we have  
$$\si_{t,r}(x,y,z)=\sum_{(\un{j},\un{c})\in \I_{\T}} \sign(\xi_{\un{j},\un{c}})\;
\si_{j_1}({c_1})\cdots \si_{j_p}({c_p})|_{1\to x,\;2\to y,\;3\to z},%
$$
where $p=\#\un{j}=\#\un{c}$. 

 
Using Lemma~4 from~\cite{Lopatin_bplp}, we can reformulate Theorem~3 from~\cite{Lopatin_bplp} for $(\T,(X,Y,Z))$ as follows: if $n=t+2r$, then
\begin{eq}\label{eq1_si} 
\frac{1}{t!\,(r!)^2 }\sum_{\xi\in S_n} \sign(\xi) 
\prod_{a\in \T^{\xi}_{\rm cl}}\tr(X_{\varphi(a)})=
\sum_{(\un{j},\un{c})\in \I_{\T}} \sign(\xi_{\un{j},\un{c}})
\;\si_{j_1}(X_{c_1})\cdots \si_{j_p}(X_{c_p}),
\end{eq}
where $p=\#\un{j}=\#\un{c}$. Consider the result of formal substitution $X\to 1$, $Y\to 2$, and $Z\to 3$ in~\Ref{eq1_si}, where $1,2,3\in\M^{\infty}$ are letters: 
\begin{eq}\label{eq01_si} 
\frac{1}{t!\,(r!)^2 }\sum_{\xi\in S_n} \sign(\xi) 
\prod_{a\in \T^{\xi}_{\rm cl}}\tr(\varphi(a))=
\sum_{(\un{j},\un{c})\in \I_{\T}} \sign(\xi_{\un{j},\un{c}})
\;\si_{j_1}(c_1)\cdots \si_{j_p}(c_p).
\end{eq}
We claim that~\Ref{eq01_si} is valid equality of two elements from $\N_{\si}^{\infty}$. Theorem~3 was proven in Sections~5, 6, and 7 of~\cite{Lopatin_bplp}. We repeat this proof without using Section~5, i.e.,  we do not apply Lemma~4 from Section~5 in the reasoning several lines before Formula~(14) from~\cite{Lopatin_bplp}. Since $\T$ has two columns, we obtain the claim. 
Lemma~\ref{lemma3_appendix} concludes the proof. 
\end{proof}

\begin{lemma}\label{lemma2_si} 
Assume that $\FF=\QQ$. Then $\si_{\un{t},\un{r},\un{s}}(x_1,\ldots,x_u,y_1,\ldots y_v,z_1,\ldots,z_w)= $
$$\frac{1}{\un{t}!\,\un{r}!\,\un{s}!}\,\,
\silin{t,r}(\underbrace{x_1,\ldots,x_1}_{t_1},\ldots,\underbrace{x_u,\ldots,x_u}_{t_u}, 
    \underbrace{y_1,\ldots,y_1}_{r_1},\ldots,\underbrace{y_v,\ldots,y_v}_{r_v},
    \underbrace{z_1,\ldots,z_1}_{s_1},\ldots,\underbrace{z_w,\ldots,z_w}_{s_w}),$$
where $t=t_1+\cdots+t_u$ and $r=r_1+\cdots+r_v=s_1+\cdots+s_w$.
\end{lemma}
\begin{proof} The left hand side of the required formula is equal to the homogeneous component of $\si_{t,r}(x_1+\cdots+x_u,y_1+\cdots+y_v,z_1+\cdots+z_w)$ of multidegree $(\un{t},\un{r},\un{s})$. Applying Lemma~\ref{lemma2a_si} to $\si_{t,r}$ and using linearity of $\silin{t,r}$, we complete the proof.   
\end{proof}

\begin{lemma}\label{lemma3_si}
Assume that $X_1,\ldots,X_t$, $Y_1,\ldots,Y_r$, $Z_1,\ldots,Z_r$ are $n\times n$ matrices and there exists a $j$ such that $Y_j$ or $Z_j$ is a symmetric matrix. Then
$$\silin{t,r}(X_1,\ldots,X_t,Y_1,\ldots,Y_r,Z_1,\ldots,Z_r)=0.$$
\end{lemma}
\begin{proof} Assume that $Y_j$ is a symmetric matrix. 
Consider a primitive closed path $\al=\be y_j$ in $\Q_{1^t,1^r,1^r}$ such that $\deg_{\ga}\al+\deg_{\ga^T}\al\leq1$ for any arrow $\ga$ from $\Q_{1^t,1^r,1^r}$. Then $\al y_j^T$ is also a primitive closed path and $\al y_j^T\not\sim \al y_j$. Lemma~\ref{lemma1_si} concludes the proof.
\end{proof}

Let us recall that $E_n$ stands for the identity $n\times n$ matrix. 

\begin{lemma}\label{lemma4_si}
Let $X_1,\ldots,X_{t-1},Y_1,\ldots,Y_r,Z_1,\ldots,Z_r$ be $n\times n$ matrices. Then 
$$\silin{t,r}(X_1,\ldots,X_{t-1},E_n,Y_1,\ldots,Y_r,Z_1,\ldots,Z_r)= a\,\silin{t-1,r}(X_1,\ldots,X_{t-1},Y_1,\ldots,Y_r,Z_1,\ldots,Z_r),$$
where $a=n-(t+2r)+1$.
\end{lemma}
\begin{proof} Denote $X_t=E_n$ and $\mu=\silin{t,r}(X_1,\ldots,X_t,Y_1,\ldots,Y_r,Z_1,\ldots,Z_r)$. Let $x_1,\ldots,z_r$ be arrows of $\Q_{1^t,1^r,1^r}$. Denote by $\psi$ the substitution $x_i\to X_i$, $y_j\to Y_j$, and $z_j\to Z_j$ for all $1\leq i\leq t$ and $1\leq j\leq r$. By Corollary~\ref{cor_si}, we have
$$\mu=\sum_{\un{\al}\in {\I}^{\rm lin}_{t,r}} %
(-1)^{\eta_{\un{\al}}} \, \psi(\tr(\un{\al}))=\mu_1+\mu_2,$$
where $\mu_1$ is the sum of all $(-1)^{\eta_{\un{\al}}} \tr(\un{\al})$ such that $\al_j\sim x_t$ for some $j$. Equalities 
$$\mu_1=\tr(E_n) \sum_{\un{\be}\in {\I}^{\rm lin}_{t-1,r}} (-1)^{\eta_{\un{\be}}} \, \psi(\tr(\un{\be}))\;\;\; \text{ and}$$
$$\mu_2=-(t+2r-1) \sum_{\un{\be}\in {\I}^{\rm lin}_{t-1,r}} (-1)^{\eta_{\un{\be}}} \, \psi(\tr(\un{\be}))$$
conclude the proof.
\end{proof}

The following remark is trivial.

\remark\label{remark_reduction_to_char_0}
Let $p>0$ be the characteristic of $\FF$ and $\FF[x_1,\ldots,x_m]$ be a polynomial ring. Define the ring homomorphism $\pi:\ZZ\to\FF$ by $\pi(1)=1_{\FF}$, where
$1_{\FF}$ stands for the unity of $\FF$. Consider
$f=\sum a_i f_i\in\FF[x_1,\ldots,x_m]$, where
$a_i\in\pi(\ZZ)$, and $f_i$ is a monomial in $x_1,\ldots,x_m$ for all $i$.
Take $b_i\in\{0,1,\ldots,p-1\}\subset\ZZ$ such that
$\pi(b_i)=a_i$ and set $h=\sum_i b_i
f_i\in\QQ[x_1,\ldots,x_m]$. Then $h=0$ implies $f=0$.

\begin{lemma}\label{lemma5_si}
Let $X,Y,Z$ be $n\times n$ matrices and $n\geq t+2r$. Then 
$$\si_{t,r}(X+\la E_n,Y,Z)=\sum\limits_{i=0}^t \binom{n-(t+2r)+i}{i} \la^i \, \si_{t-i,r}(X,Y,Z)$$ 
for any $\la\in \algA$. 
\end{lemma}
\begin{proof} By Remark~\ref{remark_reduction_to_char_0}, without loss of generality we can assume $\FF=\QQ$. Then Lemma~\ref{lemma2_si} and linearity of $\silin{t,r}$ imply 
$$\si_{t,r}(X+\la E_n,Y,Z)=\frac{1}{t!\,(r!)^2} \sum_{i=0}^t \binom{t}{i} \la^i \; \silin{t,r}(\underbrace{X,\ldots,X}_{t-i},\underbrace{E_n,\ldots,E_n}_i,Y,\ldots,Y,Z,\ldots,Z).$$
Lemma~\ref{lemma4_si} concludes the proof.
\end{proof}

\section{Proof of Theorem~\ref{theo_relations}}\label{section_proof}  

Without loss of generality we can assume that $\FF$ is an algebraically closed field (cf.~Remark~2.1 from~\cite{Zubkov_Fund_Math_01} and Remark~3.2 from~\cite{Zubkov99}).
  
Assume that $X_1=X_{n,1},\ldots,X_d=X_{n,d},Y=Y_{n},Z=Z_n$ are $n\times n$ generic matrices and entries of these matrices are $x_{ij}(1),\ldots,x_{ij}(d)$, $y_{ij}$, $z_{ij}$ $(1\leq i,j\leq n)$, respectively. Consider the mixed quiver $\G$
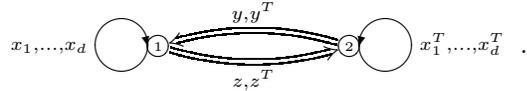
$$
\loopR{0}{0}{x_1,\ldots,x_d} %
\xymatrix@C=1cm@R=1cm{ %
\vtx{1}\ar@2@/^/@{<-}[rr]^{y,y^T} &&\vtx{2}\ar@2@/^/@{<-}[ll]^{z,z^T}\\
}%
\loopL{0}{0}{x_1^T,\ldots,x_d^T}\qquad\qquad.
$$ %
where there are $d$ loops in each of its vertices and $1^T=2$. 

Given $\al\in\M(\G)$ we define $X_{n,\al}$ in the same manner as in Section~\ref{section_intro}, using convention that
$$X_{n,\al}=
\left\{
\begin{array}{rl}
X_{n,i},&\text{if } \al=x_i\\
Y_{n},&\text{if } \al=y\\
Z_{n},&\text{if } \al=z\\
\end{array}
\right..
$$

Denote by $C_n=\FF[x_{ij}(k),y_{ij},z_{ij}\,|\,1\leq i,j\leq n,\,1\leq k\leq d]$ 
the coordinate ring of the space of mixed representations of $\G$ of dimension vector $(n,n)$. Then its algebra of {\it invariants} $J_n\subset C_n$ is generated by $\si_t(X_{n,\al})$, where $1\leq t\leq n$ and $\al$ is a primitive closed path in $\G$ (see~\cite{ZubkovI}). We define the inclusion $R_n\subset C_n$ in the natural way. 
Consider a surjective homomorphism 
$$\Upsilon_n : \N_{\si}(\G)\to J_n$$ %
defined by  $\si_t(\al) \to \si_t(X_{n,\al})$, if $t\leq n$, and $\si_t(\al) \to 0$ otherwise.
Its kernel $T_{n}$ is the ideal of relations for $J_n$. Elements of $\bigcap_{i>0} T_{i}$ are called {\it free} relations for $J_n$.  
 
Let $(\al,\be,\ga)$ be a triple of paths in $\G$. Then
it is called {\it good} if $\al'=\al''=1$, $\be$
is a path from $2$ to $1$, and $\ga$ is a path from $1$ to $2$. Since Lemma~\ref{lemma1_definition} can be reformulated for $\N_{\si}(\G)$, we have   
\begin{eq}\label{eqC_proof}
\si_{t,r}(\sum_i a_{1i} {\al_{1i}},\,\sum_j a_{2j} {\al_{2j}},\,\sum_k a_{3k} {\al _{3k}})\in\N_{\si}(\G), 
\end{eq}
where $a_{1i},a_{2j},a_{3k}\in\FF$ and  
$(\al_{1i},\al_{2j},\al_{3k})$ is a good triple for all $i,j,k$. 

\begin{theo}\label{theo_ZubII}(Zubkov,~\cite{ZubkovII}) The ideal of relations $T_{n}$ for  $J_n\simeq  \N_{\si}(\G)/T_{n}$ is generated by 
\begin{enumerate}
\item[$\bullet$] elements~\Ref{eqC_proof} for $t+2r>n$;


\item[$\bullet$] free relations for $J_n$.
\end{enumerate} 
\end{theo}

\bigskip
Given $N>n$, we define two homomorphisms $\prI_{N,n},\prII_{N,n}: C_N \to C_n$  as follows:

$$\prI_{N,n}:\quad X_{N,k}\to 
\left(
\begin{array}{cc}
X_{n,k}& 0 \\
0 & 0 \\
\end{array}
\right), \;\;
Y_{N}\to 
\left(
\begin{array}{cc}
Y_{n}& 0 \\
0 & 0 \\
\end{array}
\right), \;\; 
Z_{N}\to 
\left(
\begin{array}{cc}
Z_{n}& 0 \\
0 & 0 \\
\end{array}
\right), \;\; \text{and }
$$ %
$$
\prII_{N,n}:\quad X_{N,k}\to
\left(
\begin{array}{cc}
X_{n,k}& 0 \\
0 & 0 \\
\end{array}
\right), \;\;
Y_{N}\to 
\left(
\begin{array}{cc}
Y_{n}& 0 \\
0 & E_{N-n} \\
\end{array}
\right), \;\; 
Z_{N}\to 
\left(
\begin{array}{cc}
Z_{n}& 0 \\
0 & E_{N-n} \\
\end{array}
\right),
$$ 
where $1\leq k\leq d$ and $E_{N-n}$ stands for the identity matrix. In other words, we assume that
$$\prI_{N,n}(x_{ij}(k))=\left\{
\begin{array}{rl}
x_{ij}(k),& \text{ if } 1\leq i,j\leq n\\
0,& \text{ otherwise }\\
\end{array}
\right.$$
and so on. For short, we use the following notations: $\prI=\prI_{N,n}$ and $\prII=\prII_{N,n}$. 

Consider a mapping $\Phi_N: J_N\to J_N$ defined by:
$$\Phi_N: X_{N,k}\to X_{N,k}, Y_{N}\to E_N, Z_{N}\to E_N$$
for all $1\leq k\leq d$. 
The following diagram is commutative: 
$$ 
\begin{picture}(0,80)
\put(-37,55){%
\put(0,0){\vector(0,-1){30}}%
\put(75,0){\vector(0,-1){30}}%
\put(15,10){\vector(1,0){45}}%
\put(15,-40){\vector(1,0){45}}%
\put(135,10){\vector(-1,0){45}}%
\put(-60,10){\vector(1,0){45}}%
\put(-60,0){\vector(3,-2){45}}%
\put(135,0){\vector(-3,-2){45}}%
\put(-5,7){$J_N$}%
\put(-5,-43){$J_n$}%
\put(62,7){$R^{O(N)}$}%
\put(62,-43){$R^{O(n)}$}%
\put(-95,7){$\N_{\si}(\G)$}%
\put(145,7){$\N_{\si}$}%
\put(32,13){$\scriptstyle \Phi_N$}%
\put(32,-37){$\scriptstyle \Phi_n$}%
\put(-43,13){$\scriptstyle \Upsilon_N$}%
\put(107,13){$\scriptstyle \Psi_N$}%
\put(3,-15){$\scriptstyle \prII$}%
\put(78,-15){$\scriptstyle \prII=\prI$}%
\put(-48,-20){$\scriptstyle \Upsilon_n$}%
\put(118,-20){$\scriptstyle \Psi_n$}%

}%
\end{picture}$$%

The next lemma is a reformulation of Lemmas~3.3 and~3.4 from~\cite{ZubkovII}. 

\begin{lemma}\label{lemma01_proof}
If $N>s>0$, then 
\begin{enumerate}
\item[$\bullet$] the kernel of $\Psi_N:\Hc{s}{\N_{\si}}\to \Hc{s}{R^{O(N)}}$ is generated by free relations for $R^{O(n)}$ of degree $s$; 

\item[$\bullet$] the kernel of $\Upsilon_N:\Hc{s}{\N_{\si}(\G)}\to \Hc{s}{J_N}$ is generated by free relations for $J_n$ of degree $s$. 
\end{enumerate}
\end{lemma}

Lemma~\ref{lemma01_proof} implies that we can identify $\Hc{s}{\N_{\si}}$ with $\Hc{s}{R^{O(N)}}$ as well as $\Hc{s}{\N_{\si}(\G)}$ with $\Hc{s}{J_N}$ modulo free relations.
Let us formulate Proposition~3.2 from~\cite{ZubkovII}.

\begin{lemma}\label{lemma02_proof}
The ideal $K_n$ is generated by elements $\Phi_N(h_f)$, where 
\begin{enumerate}
\item[$\bullet$] $f$ ranges over elements~\Ref{eqC_proof} with $t+2r>n$; 

\item[$\bullet$] $N>\deg{f}=s$ and $N\gg n$ is big enough, so we can assume that $f\in J_N$;

\item[$\bullet$] $h=h_f\in J_N$ such that $\prII(h)=0$, $\deg{h}=s$, and $\hc{s}{h}=\hc{s}{f}$. 
\end{enumerate}
Here $\Phi_N(h_f)$ is considered as an element of $\N_{\si}$.
\end{lemma}
\bigskip

\noindent{}Note that in Lemma~\ref{lemma02_proof} $f$ can be non-homogeneous.

Given $t,r\geq0$ satisfying $t+2r>n$, we consider an element~\Ref{eqC_proof} 
$$f=\si_{t,r}(f_1,f_2,f_3),$$ 
where $f_1 = \sum_i a_{1i} {\al_{1i}}$, $f_2=\sum_j a_{2j} {\al_{2j}}$, and $f_3=\sum_k a_{3k} {\al_{3k}}$.  There exists an $N\gg n$ such that $N>t+2r$, $N>t+n$, and we can assume that $f\in J_N$. In what follows we will write $f_1$ instead of $\sum_i a_{1i} X_{N,\al_{1i}}$ and similarly for $f_2$ and $f_3$.

By Theorem~\ref{theo_ZubII}, we have 
\begin{eq}\label{eq3_proof}
\prI(f)=0.
\end{eq}
But we can not claim that $\prII(f)=0$. We will construct $h=h_f\in C_N$ that is ``close'' to $f$ and $\prII(h)=0$ (see Lemma~\ref{lemma2_proof} below). 

We can rewrite $f_1$ as $f_1=f_{11}+f_{12}$, where $f_{11}$ is a sum of all $a_{1i} X_{N,\al_{1i}}$ such that $\al_{1i}$ contains an arrow $x_k$ or $x_k^T$ for some $k$. Similarly, we can rewrite $f_2$, $f_3$ as $f_2=f_{21}+f_{22}$ and $f_3=f_{31}+f_{32}$, where $f_{21}$ and $f_{31}$ contain all summands with $x_k$ or $x_k^T$. Let $\la_l\in\FF$ be the sum of coefficients of summands of $f_{l,2}$ for $l=1,2,3$. Note that
\begin{eq}\label{eq2_proof}
\prII(f_l)=\prI(f_l)+\la_l E_{N,n}
\end{eq} 
for $l=1,2,3$. Here $E_{N,n}$ stands for the following $N\times N$ matrix
$$
\left(
\begin{array}{cc}
0& 0 \\
0 & E_{N-n} \\
\end{array}
\right). 
$$ 


\begin{lemma}\label{lemma1_proof}
$$\prII(f)=\sum_{i=0}^t \binom{N-n}{i} \prI(\si_{t-i,r}(f_1,f_2,f_3))\; 
\la_1^i.
$$
\end{lemma}
\begin{proof} By Remark~\ref{remark_reduction_to_char_0}, without loss of generality we can assume $\FF=\QQ$. By~\Ref{eq2_proof}, we have
$$\prII(f)=\si_{t,r}(\,\prI(f_1)+\la_1 E_{N,n},\;
\prI(f_2)+ \la_2 E_{N,n},\; \prI(f_3)+ \la_3 E_{N,n}\,).$$ %
Lemmas~\ref{lemma2a_si},~\ref{lemma3_si} together with linearity of $\silin{t,r}$ imply $\prII(f)=$
$$\frac{1}{t!\,(r!)^2}\,\, \sum_{i=0}^t \binom{t}{i}
\silin{t,r}(\, \underbrace{\prI(f_1),\ldots,\prI(f_1)}_{t-i},\, \underbrace{E_{N,n},\ldots,E_{N,n}}_i,\, 
\prI(f_2),\ldots,\prI(f_2),\, \prI(f_3),\ldots,\prI(f_3))\, \la_1^i.$$
Since $\prI(f_l) E_{N,n}=0$ for all $l$, applying Lemma~\ref{lemma1_si} we obtain
$$\prII(f)=\frac{1}{(r!)^2}\,\, \sum_{i=0}^t \frac{1}{i!\,(t-i)!}\,\,
\silin{t-i,r}(\prI(f_1),\ldots,\prI(f_1), \prI(f_2),\ldots,\prI(f_2), \prI(f_3),\ldots,\prI(f_3))\,\la_1^i a,
$$ 
where $a=\silin{i,0}(E_{N,n},\ldots,E_{N,n})$. By Lemma~\ref{lemma2a_si} and Remark~\ref{rem01}, we have 
$$a=i!\,\si_{i,0}(E_{N,n},0,0)=i!\,\si_i(E_{N,n})=i!\,\binom{N-n}{i}.$$
Lemma~\ref{lemma2a_si} concludes the proof.
\end{proof}

\begin{lemma}\label{lemma2_proof}
There exist $b_0,b_1,\ldots,b_t\in\ZZ$ such that 
$$h=h_f=\sum_{i=0}^t b_i\la_1^i\si_{t-i,r}(f_1,f_2,f_3)\in J_N$$
satisfies the following conditions:
\begin{enumerate}
 \item[a)] $\prII(h)=0$; 
 \item[b)] $\deg h=s$ and $\hc{s}{h}=\hc{s}{f}$, where $s=\deg f$.
\end{enumerate}
\end{lemma}
\begin{proof} We set $b_0=1$. By Lemma~\ref{lemma1_proof}, we have %
$$\prII(h)=\sum_{i=0}^t b_i\,\la_1^i \sum_{j=0}^{t-i} \binom{N-n}{j} \prI(\si_{t-i-j,r}(f_1,f_2,f_3))\, \la_1^j.$$
We substitute $k$ for $t-i-j$ and obtain
$$\prII(h)=\sum_{k=0}^t \prI(\si_{k,r}(f_1,f_2,f_3))\,\la_1^{t-k}\,\sum_{i=0}^{t-k} b_i\binom{N-n}{t-i-k}.$$
The following system of linear equations 
\begin{eq}\label{eq5_proof}
\sum_{i=0}^{t-k} b_i\binom{N-n}{t-i-k}=0, \text{ where } 0\leq k<t
\end{eq}
with respect to $b_1,\ldots,b_t$ is triangular and has $1$ on the main diagonal. Hence this system has a solution $b_1,\ldots,b_t\in\ZZ$.  Using Formula~\Ref{eq3_proof}, we can see that the equality $\prII(h)=0$ holds. Obviously, $\deg h=s$ and $\hc{s}{h}=\hc{s}{b_0 \si_{t,r}(f_1,f_2,f_3)}=\hc{s}{f}$.
\end{proof}

\begin{lemma}\label{lemma3_proof}
Assume that $b_0=1$ and $b_1,\ldots,b_t$ satisfy system of linear equations~\Ref{eq5_proof}. Then 
$$\sum_{i=0}^{t-k} b_{t-i-k} \binom{N-(l+2r)}{i}=0$$
for $0\leq k\leq l \leq n-2r$.
\end{lemma}
\begin{proof} Denote the left hand side of the required formula by $b_{k,l}$.

The proof is by decreasing induction on $l$. If $l=n-2r$, then $b_{k,l}=0$ by~\Ref{eq5_proof}. 

Let $l<n-2r$. Then formula
$$\binom{q-1}{p}+\binom{q-1}{p-1}=\binom{q}{p},
$$
where $0\leq p<q$, implies
$$b_{k,l}=b_{k,l+1} + \sum_{i=0}^{t-k} b_{t-i-k} \binom{N-(l+2r+1)}{i-1}.$$
We substitute $j$ for $i-1$ and obtain $b_{k,l}=b_{k,l+1}+b_{k+1,l+1}$. Induction hypothesis concludes the proof.
\end{proof}

Denote by $I_N$ the ideal of $R^{O(N)}$ generated by $\si_{p,q}(\sum_i c_{1i} X(\ga_{1i}),\,\sum_j c_{2j} X(\ga_{2j}),\,\sum_k c_{3k} X(\ga_{3k}))$, where $X(\ga)$ stands for $X_{N,\ga}$, $p+2q>n$, $\ga_{1i},\ga_{2j},\ga_{3k}$ are monomials in $x_1,\ldots,x_d,x_1^T,\ldots,x_d^T$, and $c_{1i},c_{2j},c_{3k}\in \FF$. 

\begin{lemma}\label{lemma4_proof}
Assume that $h\in J_N$ is the element from Lemma~\ref{lemma2_proof}. Then $\Phi_N(h)$ belongs to $I_N$.
\end{lemma}
\begin{proof} For short, denote $\Phi_N(f_{l1})=g_{l}$ for $l=1,2,3$. We have $\Phi_N(f_{l})=g_{l}+\la_l E_N$ for all $l$. 

To begin with, we assume that $\FF=\QQ$. Using Lemmas~\ref{lemma2a_si} and~\ref{lemma3_si} together with linearity of $\silin{t,r}$ in the same manner as in the proof of Lemma~\ref{lemma1_proof}, we obtain
$$\Phi_N(h)=\sum_{i=0}^t b_i\, \la_1^i\, \si_{t-i,r}(g_1+\la_1 E_N, g_2, g_3).$$
Lemma~\ref{lemma5_si} implies
\begin{eq}\label{eq6_proof}
\Phi_N(h)=\sum_{k=0}^t \si_{k,r}(g_1,g_2,g_3)\, \la_1^{t-k} \sum_{i=0}^{t-k} b_{t-i-k} \binom{N-(k+2r)}{i}.
\end{eq}
In the general case Remark~\ref{remark_reduction_to_char_0} shows that Formula~\Ref{eq6_proof} is also valid.

If $k+2r\leq n$, then 
\begin{eq}\label{eq4_proof}
\sum_{i=0}^{t-k} b_{t-i-k} \binom{N-(k+2r)}{i}=0
\end{eq}
by Lemma~\ref{lemma3_proof}. Thus $\Phi_N(h)$ belongs to $I_N$.
\end{proof}

\begin{lemma}\label{lemma5_proof}
Any element $g=\si_{t,r}(\sum_i c_{1i} {\ga_{1i}},\,\sum_j c_{2j} {\ga_{2j}},\,\sum_k c_{3k} {\ga_{3k}})$ of $\N_{\si}$, where $t+2r>n$ and $c_{1i},c_{2j},c_{3k}\in\FF$, 
$\,\ga_{1i},\ga_{2j},\ga_{3k}\in\M$, belongs to $K_n$.
\end{lemma}
\begin{proof}
There exists a $g'\in\N_{\si}(\G)$ such that $\Phi_N(\Upsilon_N(g'))=\Psi_N(g)$. Since  $\Upsilon_n(g')=0$ by Theorem~\ref{theo_ZubII}, we obtain 
$$\Psi_n(g)=p\circ \Psi_N(g)=p\circ \Phi_N \circ \Upsilon_N (g')=\Phi_n\circ \Upsilon_n(g')=0.$$
\end{proof}

Lemmas~\ref{lemma02_proof},~\ref{lemma4_proof}, and~\ref{lemma5_proof} conclude the proof of Theorem~\ref{theo_relations}.

\section{Free relations}\label{section_free}

In this section it is convenient to write $x_1,\ldots,x_d$, $x_1^T,\ldots,x_d^T$ for letters $1,\ldots,d$, $1^T,\ldots,d^T$ that are free generators of $\M$.
Similarly to Section~\ref{section_intro}, we denote by $\M^{\infty}$ the monoid freely generated by letters $x_1,x_2,\ldots$, $x_1^T,x_2^T,\ldots$ and define $\N^{\infty}$ ($\N_{\si}^{\infty}$, respectively) similarly to $\N$ ($\N_{\si}$, respectively).  

\definition{of $\Lin(f)$} Given $\NN^d$-homogeneous $f\in\N_{\si}$ with $\mdeg{f}=\un{t}=(t_1,\ldots,t_d)$, we define its {\it complete linearization}  $\Lin(f)\in\N_{\si}^{\infty}$ as follows. Let $h\in\N_{\si}^{\infty}$ be the result of substitution $x_i\to \sum_{j=0}^{{t_i} - 1}a_{ij} x_{i+jd}$ ($1\leq i\leq d$) in $f$, where $a_{ij}\in\FF$, and we consider $h$ as a polynomial in $a_{ij}$. We denote the coefficient of $a_{11}\cdots a_{1 t_1}\cdots a_{d 1}\cdots a_{d t_d}$ in $h$ by $\Lin(f)$.

\definition{of $\Glue(\al)$} Given $\al\in\M^{\infty}$, we define $\Glue(\al)\in\M$ as the result of substitution $x_{i+jd}\to x_i$ ($1\leq i\leq d$ and $j>0$) in $\al$.

\example\label{ex1_free} If $f=\si_2(x_1)$, then~\Ref{eq_si2} implies that $\Lin(f)=-\tr(x_1 x_{d+1})+\tr(x_1)\tr(x_{d+1})$.

If $f=\tr(x_1)^3$, then $\Lin(f)=6\tr(x_1)\tr(x_{d+1})\tr(x_{2d+1})$.

\bigskip
Given $f=\si_{j_1}(\al_1)\cdots\si_{j_p}(\al_p)\in\N_{\si}$, where $\al_1,\ldots,\al_p\in\N$, we define $c_f$ and $e_f$ as follows. We consider a subset $S$ of the set of pairs $(j_i,\al_i)$, $1\leq i\leq p$, such that elements of $S$ are pairwise different with respect to $\sim$, where $(j_i,\al_i)\sim (j_k,\al_k)$ if and only if $j_i=j_k$ and $\al_i\sim\al_k$.  Given $(j,\al)\in S$, we denote $c_{(j,\al)}=\#\{(j_i,\al_i)\,|\,j_i=j,\al_i\sim\al, 1\leq i\leq p\}$. We set  $c_f=\prod_{(j,\al)\in S} c_{(j,\al)}!$ and $e_f=p$. Denote by $\ZZ$ the subring of $\FF$
generated by $\pm1$.

The proof of the following lemma is similar to the proof of Lemma~\ref{lemma1_si}. 

\begin{lemma}\label{lemma2_free}
Let $f=\si_{j_1}(\al_1)\cdots\si_{j_p}(\al_p)\in\N_{\si}$, where $\al_1,\ldots,\al_p\in\N$, and $\un{t}=\mdeg{f}$.  Then 
$$\Lin(f)=\pm 2^r c_f\, \sum \tr(\ga_1)\cdots\tr(\ga_p)+c_f h,$$
where 
\begin{enumerate}
\item[$\bullet$] $\ga_1,\ldots,\ga_p\in\N^{\infty}$ are representatives of pairwise different $\sim$-equivalence classes and $r\geq0$;

\item[$\bullet$] we have $\deg_{x_l}(\ga_1\cdots\ga_p)+\deg_{x_l^T}(\ga_1\cdots\ga_p)=\left\{
\begin{array}{cc}
1,&\text{ if }l=i+jd \text{ for } 1\leq i\leq d \text{ and } 0\leq j<t_i\\
0,&\text{ otherwise}\\
\end{array}
\right.$ for all $l>0$;

\item[$\bullet$] $\Glue(\ga_i)=\al_i^{j_i}$ for all $1\leq i\leq p$.

\item[$\bullet$] $h=\sum_{k} \be_k h_k$, where $\be_k\in\ZZ$ and $h_k$ is a product of $p+1$ or more traces.  
\end{enumerate}
\end{lemma}

\begin{lemma}\label{lemma1_free}
If $\Char{\FF}>0$, then the ideal of free relations lies in the ideal of $\N_{\si}$ generated by $\si_t(\al)^{\Char{\FF}}$, where $\al$ ranges over $\N$ and $t>0$. If $\Char{\FF}=0$, then the only free relation is trivial. 
\end{lemma} 
\begin{proof} 
Consider an $\NN^d$-homogeneous free relation $f\in\N_{\si}$.

\medskip
{\bf 1.} Let $\mdeg{f}=\un{t}$ satisfy $t_i\leq1$ for all $1\leq i\leq d$. We assume that $f\neq0$. Then $f=\sum_k a_k f_k$, where $a_k\in\FF$, $a_k\neq0$, and $f_k$ are pairwise different products of traces.  Denote $n=\deg{f}$ and denote by $e_{i,j}$ the $n\times n$ matrix whose $(i,j)^{\rm th}$ entry is $1$ and any other entry is $0$. Let $f_1=\tr(\al_1)\cdots\tr(\al_p)$ for some $\al_1,\ldots,\al_p\in\N$. Given $\al_1=y_1\cdots y_q$, where $y_1,\ldots,y_q$ are letters, we set $Y_i=e_{i,i+1}$ for $1\leq i<q$ and $Y_q=e_{q,1}$. Considering $\al_2,\ldots,\al_p$, we define $Y_i$ for $q<i\leq \deg{f_1}$. 
Since $f$ is a free relation, the substitution $y_i\to Y_i$ ($1\leq i\leq\deg{f_1}$) implies that $a_1=0$; a contradiction. Thus, $f=0$.

\medskip
{\bf 2.} Assume that $f$ is non-zero. Since $\Lin(f)$ is a free relation, part~1 implies that $\Lin(f)=0$. We have $f=\sum_i a_i f_i+\sum_j b_j h_j$, where $a_i,b_j\in\FF$, $a_i,b_j\neq0$, $f_i,h_j$ are pairwise different products of $\si_t$ ($t>0$) such that for some $q$
\begin{enumerate}
\item[$\bullet$] $e_{f_i}=q$ for all $i$;

\item[$\bullet$] $e_{h_j}>q$ for all $j$. 
\end{enumerate}
Linearity of $\Lin$ and Lemma~\ref{lemma2_free} imply that $c_{f_i}=0$. Hence we obtain that $f_i$ lies in the required ideal of $\N_{\si}$ in case $\Char{\FF}>0$; and we obtain a contradiction in case $\Char{\FF}=0$. Since $\Lin(f_i)=0$, we have $\Lin(h)=0$ for $h=\sum_j b_j h_j$. Then we repeat the above reasoning for $h$ and so on.
\end{proof}

Given an $\NN$-graded algebra $\algA$, denote by $\algA^{+}$ the subalgebra generated by elements of $\algA$ of positive degree. It is easy to see that a set $\{a_i\} \subseteq \algA$ is a minimal (by inclusion) homogeneous system of generators (m.h.s.g.) for $\algA$~if and only if $\{\ov{a_i}\}$ is a basis for $\ov{\algA}={\algA}/{(\algA^{+})^2}$ and $\{a_i\}$ are homogeneous. Let us recall that an element $a\in\algA$ is called {\it
decomposable} if it belongs to the ideal $(\algA^{+})^2$. 
Therefore the
least upper bound for the degrees of elements of a m.h.s.g.~for $R^{O(n)}$ is equal to
the maximal degree of indecomposable invariants and we denote it by $D_{\rm max}$. Theorem~\ref{theo_relations} together with Lemma~\ref{lemma1_free} imply the following corollary.

\begin{cor}\label{cor_free} If $\Char{\FF}\neq2$, then the ideal of relations $\ov{K_{n}}$ of $\ov{R^{O(n)}}\simeq  \ov{\N_{\si}}/\ov{K_{n}}$ is generated by $\ov{\si_{t,r}(\al,\be,\ga)}$, where $t+2r>n$ and 
$\al,\be,\ga$ range over $\M_{\FF}$.
\end{cor}

\bigskip
\noindent As an application of Corollary~\ref{cor_free} we obtained the following result in~\cite{Lopatin_O3}. 

\begin{theo}\label{theo_free} Let $n=3$ and $d\geq1$. Then
\begin{enumerate} 
\item[$\bullet$] If $\Char\FF =3$, then $2d+4\leq D_{\rm max}\leq 2d+7$.

\item[$\bullet$] If $\Char\FF \neq 2,3$, then $D_{\rm max}=6$.
\end{enumerate}
\end{theo}

As about matrix $GL(n)$-invariants in case $n=3$, its minimal system of generators was explicitly calculated in~\cite{Lopatin_Comm1} and~\cite{Lopatin_Comm2}.

\section{Appendix: other definitions of $\si_{t,r}$}\label{section_appendix}

In this section we assume that $\algA$ is a commutative unitary algebra over the field $\FF$ and all matrices are considered over $\algA$. In what follows we recall some definitions from~\cite{Lopatin_bplp}. Note that in this section we consider only rectangular tableaux with two columns whereas in~\cite{Lopatin_bplp} tableaux with arbitrary number of columns of any length were defined.   

\begin{definition}{of a tableau with substitution} 
A pair $(\T,(X_1,\ldots,X_d))$ is called a {\it tableau with substitution} of dimension $(n,n)$, if 
\begin{enumerate}
\item[$\bullet$] $\T$ is a rectangular tableau with two columns and $n$ rows.  The tableau $\T$ is filled with arrows in such a way that an {\it arrow} goes from one cell of the tableau into another one, and each cell of the tableau is either the head or the tail of one and only one arrow. We write $a\in \T$ for an  arrow $a$ from $\T$. Given an arrow $a\in \T$, denote by $a'$ and $a''$ the columns containing the head and the tail of $a$, respectively. Similarly, denote by $'a$ the row containing the head of $a$, and denote by $''a$ the row containing the tail of $a$. Schematically this is depicted as 
$$
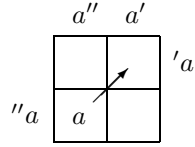
\begin{picture}(100,50)
\put(50,30){%
\put(0,0){\rectangle{10}{10}\put(-3,-3){}}%
\put(20,0){\rectangle{10}{10}\put(-3,-3){}}%
\put(0,-20){\rectangle{10}{10}\put(-3,-3){$a$}\put(5,5){\vector(1,1){13}}}%
\put(20,-20){\rectangle{10}{10}\put(-3,-3){}}%
\put(-3,15){$a''$}%
\put(17,15){$a'$}%
\put(35,-3){$'a$}%
\put(-26,-23){$''a$}%
}%
\end{picture}$$

\item[$\bullet$] $\varphi$ is a fixed mapping from the set of arrows of $\T$ onto $\{1,\ldots,d\}$ that satisfies the following property:

\begin{enumerate}
\item[] if $a,b\in \T$ and $\varphi(a)=\varphi(b)$, then $a'=b'$, $a''=b''$;
\end{enumerate}

\item[$\bullet$] $X_1,\ldots,X_d$ are $n\times n$ matrices.
\end{enumerate}
\end{definition} 

\begin{example}\label{ex1_appendix}
Let $X,Y,Z$ be $(t+2r)\times (t+2r)$ matrices and let $\T=\T_{t,r}$ be the tableau:
$$
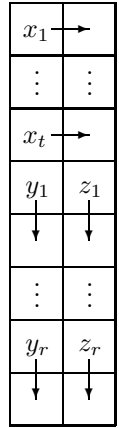
\begin{picture}(100,160)
\put(50,150){%
\put(0,0){\rectangle{10}{10}\put(-5,-3){$x_1$}\put(6,0){\vector(1,0){14}}}%
\put(20,0){\rectangle{10}{10}}%
\put(0,-20){\rectangle{10}{10}\put(-1,-5){$\vdots$}}%
\put(20,-20){\rectangle{10}{10}\put(-1,-5){$\vdots$}}%
\put(0,-40){\rectangle{10}{10}\put(-5,-3){$x_t$}\put(6,0){\vector(1,0){14}}}%
\put(20,-40){\rectangle{10}{10}}%

\put(0,-60){\rectangle{10}{10}\put(-4,-2){$y_1$}\put(0,-5){\vector(0,-1){15}}}%
\put(0,-80){\rectangle{10}{10}}%
\put(0,-100){\rectangle{10}{10}\put(-1,-5){$\vdots$}}%
\put(0,-120){\rectangle{10}{10}\put(-4,-2){$y_r$}\put(0,-5){\vector(0,-1){15}}}%
\put(0,-140){\rectangle{10}{10}}%

\put(20,-60){\rectangle{10}{10}\put(-4,-2){$z_1$}\put(0,-5){\vector(0,-1){15}}}%
\put(20,-80){\rectangle{10}{10}}%
\put(20,-100){\rectangle{10}{10}\put(-1,-5){$\vdots$}}%
\put(20,-120){\rectangle{10}{10}\put(-4,-2){$z_r$}\put(0,-5){\vector(0,-1){15}}}%
\put(20,-140){\rectangle{10}{10}}%
\put(40,-70){.}
}%
\end{picture}
$$ %

\noindent We define $\varphi$ as follows: $\varphi(x_i)=1$, $\varphi(y_j)=2$, $\varphi(z_j)=3$ for $1\leq i\leq t$ and $1\leq j\leq r$. Then $(\T,(X,Y,Z))$ is a tableau with substitution of dimension $(t+2r,t+2r)$.%
\end{example}

\begin{example}\label{ex1b_appendix}
Let $X_1,\ldots,X_t$, $Y_1,\ldots,Y_r$, $Z_1\ldots,Z_r$ be $(t+2r)\times (t+2r)$ matrices and let $\T$ be the tableau from Example~\ref{ex1_appendix}. We define $\varphi$ as follows: $\varphi(x_i)=i$, $\varphi(y_j)=t+j$, $\varphi(z_j)=t+r+j$ for $1\leq i\leq t$ and $1\leq j\leq r$. Then $(\T,(X_1,\ldots,X_t,Y_1,\ldots,Y_r,Z_1\ldots,Z_r))$ is a tableau with substitution of dimension $(t+2r,t+2r)$.%
\end{example}

\begin{definition}{of $\bpf_{\T}(X_1,\ldots,X_d)$} Let $X_1,\ldots,X_d$ be $n\times n$ matrices and let $(\T,(X_1,\ldots,X_d))$ be a tableau with substitution of dimension $(n,n)$. Then we define 
$$\bpf_{\T}(X_1,\ldots,X_d)=\sum \sign(\pi_1) \sign(\pi_2) %
\prod_{a\in \T} (X_{\varphi(a)})_{\pi_{a''}(''a),\pi_{a'}('a)},$$ %
where $(X_{\varphi(a)})_{ij}$ stands for $(i,j)^{\rm th}$ entry of the matrix $X_{\varphi(a)}$ and the sum ranges over permutations $\pi_1,\pi_2\in S_n$ such that for any
$a,b\in \T$ the conditions $\varphi(a)=\varphi(b)$ and $''a<{}''b$
imply that $\pi_{i}({}''a)<\pi_{i}({}''b)$ for $i=a''=b''$.
For $\FF=\QQ$ there exists a more convenient formula
$$\bpf_{\T}(X_1,\ldots,X_d)=\frac{1}{t!\,(r!)^2}  \sum_{\pi_1,\pi_2\in S_n}\sign(\pi_1)\sign(\pi_2) %
\prod_{a\in \T} (X_{\varphi(a)})_{\pi_{a''}(''a),\pi_{a'}('a)}.
$$
\end{definition}

\begin{example}\label{ex1c_appendix} 
Assume that $(\T,(X,Y,Z))$ is the tableau with substitution from Example~\ref{ex1_appendix}. Then 
$$\DP_{r,r}(X,Y,Z)=\bpf_{\T}(X,Y,Z),$$ %
where the determinant-pfaffian $\DP_{r,r}(X,Y,Z)$ was introduced in~\cite{LZ1}. 
\end{example}
\bigskip

The decomposition formula was formulated in Theorem~3,~\cite{Lopatin_bplp} (see~\Ref{eq_decomp_appendix} below). Assume that $(\T,(X,Y,Z))$ is the tableau with substitution from Example~\ref{ex1_appendix}. Denote by $\dec_{\T}(X,Y,Z)$ the right hand side of the decomposition formula, applied to the tableau with substitution $(\T,(X,Y,Z))$. Given $x,y,z\in\M_{\FF}$, we assume that $\dec_{\T}(x,y,z)\in\N_{\si}$ stands for the result of formal substitution $X\to x$, $Y\to y$, and $Z\to z$ in $\dec_{\T}(X,Y,Z)$.  

\begin{lemma}\label{lemma1_appendix}%
We have $\si_{t,r}(x,y,z)=\dec_{\T}(x,y,z)$ for $x,y,z\in\M_{\FF}$, where $\T=\T_{t,r}$ is defined in Example~\ref{ex1_appendix}.

In particular,  $\si_{t,r}(X,Y,Z)=\DP_{r,r}(X,Y,Z)$ for $(t+2r)\times(t+2r)$ matrices $X,Y,Z$. 
\end{lemma}
\bigskip

To prove this lemma we need we need additional definitions from~\cite{Lopatin_bplp}.   

Let $(\T,(X_1,\ldots,X_d))$ be a tableau with substitution of dimension $(n,n)$. Consider $\M^{\infty}$ and $\N_{\si}^{\infty}$ from Section~\ref{section_free}, where letters are $1,2,\ldots$, $1^T,2^T,\ldots$ Given $a\in \T$, we consider
$\varphi(a)\in\{1,\ldots,d\}$ as an element of $\M^{\infty}$. For $u\in\M^{\infty}$ define the matrix $X_u$ in the same way as in Section~\ref{section_intro}. Here we assume that $X_{i}=X_{i^T}=0$ for $i>d$.

For an arrow $a\in \T$ denote by $a^T$ the {\it transpose arrow},
i.e., by definition $(a^T)''=a'$, $(a^T)'=a''$, $''(a^T)={}'a$,
$'(a^T)={}''a$, and $\varphi(a^T)=\varphi(a)^T\in\M^{\infty}$. We write $a\stackrel{T}{\in} \T$ if $a$ is an arrow or a transpose arrow of $\T$.

\begin{definition}{of paths} We say that
$a_1,a_2\stackrel{T}{\in}\T$ are {\it successive} in $\T$,  if $a_1'\neq a_2''$ and $'a_1={}''a_2$.

A word $a=a_1\cdots a_s$, where
$a_1,\ldots,a_s\stackrel{T}{\in}\T$, is called a {\it path} in $\T$, if $a_i$, $a_{i+1}$ are
successive for any $1\leq i\leq s-1$. In this case by definition
$\varphi(a)=\varphi(a_1) \cdots \varphi(a_s)\in \M^{\infty}$ and
$a^T=a_s^T\cdots a_1^T$ is a path in $\T$; we denote
$a_s',{}'a_s,a_1'',{}''a_1$, respectively, by $a',{}'a,a'',{}''a$,
respectively. 

A path $a_1\cdots a_s$ is {\it closed} if $a_s,a_1$ are successive. 

Consider words $a=a_1\cdots a_p$, $b=b_1\cdots b_q$, where
$a_1\cdots a_p,b_1\cdots b_q\stackrel{T}{\in}\T$ (in particular,
both $a,b$ might be paths in $\T$). We write $a\sim b$ if there is a cyclic
permutation $\pi\in S_p$ such that $a_{\pi(1)}\cdots
a_{\pi(p)}$ is equal to $b$ or $b^T$. We will use similar notations for elements of $\M^{\infty}$, for sets of paths, and for subsets of $\M^{\infty}$.

Denote by  $\T_{\rm cl}$ a set of representatives of closed paths in $\T$ with respect to the $\sim$-equivalence.
\end{definition}

\begin{example}\label{ex2_appendix} Let %
$$
\begin{picture}(100,60)%
\put(50,50){%
\put(0,0){\rectangle{10}{10}\put(-3,-3){$a$}\put(5,0){\vector(1,0){15}}}%
\put(20,0){\rectangle{10}{10}}%
\put(0,-20){\rectangle{10}{10}}%
\put(20,-20){\rectangle{10}{10}}%
\put(0,-40){\rectangle{10}{10}\put(-3,-3){$b$}\put(5,5){\vector(1,1){13}}}%
\put(20,-40){\rectangle{10}{10}\put(-3,-3){$c$}\put(-5,5){\vector(-1,1){13}}}%
}%
\end{picture}$$
be a fragment of tableau. Then $a$ and $bc^T$ are closed paths.
\end{example}

\begin{definition}{of $\T^{\tau}$} %
Given $\tau\in S_{n}$, we permute the cells of the
$2^{\rm nd}$ column of $\T$ by $\tau$ and denote the resulting tableau
by $\T^{\tau}$. The arrows of $\T^{\tau}$ are $\{a^{\tau}\,|\,a\in
\T\}$, where $\varphi(a^{\tau})=\varphi(a)$, $(a^{\tau})''=a''$,
$(a^{\tau})'=a'$, and %
$$''(a^{\tau})=\left\{
\begin{array}{cl}
''a,& \text{if } a''=1\\
\tau(''a),& \text{if } a''=2\\
\end{array}
\right. ,\quad%
'(a^{\tau})=\left\{
\begin{array}{cl}
'a,&\text{if } a'=1\\
\tau('a),&\text{if } a'=2\\
\end{array}
\right..
$$
Obviously, $(\T^{\tau},(X_1,\ldots,X_d))$ is a tableau with substitution. For short we will write $\T_{\rm cl}^{\tau}$ instead of $(\T^{\tau})_{\rm cl}$.
\end{definition}
\bigskip

We use notation $\{\ldots\}_{m}$ for {\it multisets}, i.e., given an equivalence $=$ on a set $S$ and
$a_1,\ldots,a_p,b_1,\ldots,b_q\in S$, we write $\{a_1,\ldots,a_p\}_{m}=\{
b_1,\ldots,b_q\}_{m}$ if and only if $p=q$ and
$$\#\{1\leq j\leq p\,|\,a_j=a_i\}=\#\{1\leq j\leq p\,|\,b_j=a_i\}$$ %
for any $1\leq i\leq p$.
\smallskip


\begin{definition}{of $\I_{\T}$} Let $\underline{j}\in\NN^p$ and 
$\underline{c}=(c_1,\ldots,c_p)$, where $c_1,\ldots,c_p\in\M^{\infty}$
are primitive and pairwise different with respect to $\sim$ . Then
$(\underline{j},\underline{c})$ is called a {\it $\T$-pair}.

A $\T$-pair $(\underline{j},\underline{c})$ is called {\it
$\T$-admissible} if for some
$\xi=\xi_{\underline{j},\underline{c}}\in S_{n}$ the
following equivalence of multisets holds: %
$$\varphi(\T^{\xi}_{\rm cl})\sim%
\{\underbrace{c_1,\ldots,c_1}_{j_1},\ldots,\underbrace{c_p,\ldots,c_p}_{j_p}\}_{m}.$$
We write
$(\underline{j^0},\underline{c^0})\sim(\underline{j},\underline{c})$
and say that these pairs are equivalent if and only if %
$$%
\{\underbrace{c_1^0,\ldots,c_1^0}_{j_1^0},\ldots,
\underbrace{c_p^0,\ldots,c_p^0}_{j_p^0}\}_{m}%
\sim%
\{\underbrace{c_1,\ldots,c_1}_{j_1},\ldots,\underbrace{c_p,\ldots,c_p}_{j_p}\}_{m}.
$$
If $(\underline{j^0},\underline{c^0})\sim (\underline{j},\underline{c})$ and
$(\underline{j},\underline{c})$ is $\T$-admissible, then the
pair $(\underline{j^0},\underline{c^0})$ also has the
same property, since we can take
$\xi_{\underline{j^0},\underline{c^0}}=\xi_{\underline{j},\underline{c}}$.
Denote by $\I_{\T}$ a set of representatives of $\T$-admissible
pairs with respect to the $\sim$-equivalence.
\end{definition}
\bigskip

The decomposition formula states
\begin{eq}\label{eq_decomp_appendix}
\bpf_{\T}(X_1,\ldots,X_d)=\sum_{(\un{j},\un{c})\in \I_{\T}} \sign(\xi_{\un{j},\un{c}})\;
\si_{j_1}(X_{c_1})\cdots \si_{j_p}(X_{c_p}),%
\end{eq}
where $p=\#\un{j}=\#\un{c}$. Note that $\xi_{\un{j},\un{c}}$ is not 
unique for a representative of the $\sim$-equivalence class of a
$\T$-admissible pair 
$(\un{j},\un{c})$, but $\sign(\xi_{\un{j},\un{c}})$ is unique (see Lemma~\ref{lemma2_appendix} below). Hence
\begin{eq}\label{eq_decomp2_appendix}
\dec_{\T}(x,y,z)=\sum_{(\un{j},\un{c})\in \I_{\T}} \sign(\xi_{\un{j},\un{c}})\;
\si_{j_1}({c_1})\cdots \si_{j_p}({c_p})|_{1\to x,\;2\to y,\;3\to z}%
\end{eq}
for $x,y,z\in\M_{\FF}$. 
 
\begin{proof_without_dot} {\it $\!\!\!\!$ of Lemma~\ref{lemma1_appendix}.} %
Given an arrow or transpose arrow $a$ of $\T$, we can consider $\varphi(a)\in\M^{\infty}$  as an arrow of the quiver $\Q$ from Section~\ref{section_definition} as follows: 
$\varphi(x_i)$, $\varphi(x_i^T)$ are loops in vertices $1$ and $2$ of $\Q$, respectively, $\varphi(y_j)$ goes from vertex $2$ to vertex $1$ and so on. Similarly, for any word $a$ in arrows and transpose arrows of $\T$ we can consider $\varphi(a)$ as a word in arrows of $\Q$. Note that
\begin{enumerate}
\item[$\bullet$] if $a$ is a path in $\T$, then $\varphi(a)$ is a path in $\Q$;

\item[$\bullet$] if $a\in\T_{\rm cl}$, then $\varphi(a)$ is a closed path in $\Q$.
\end{enumerate}
Hence for any $\T$-admissible pair $(\un{j},\un{c})$ we can consider $c_1,\ldots,c_p$ as closed paths in $\Q$. Therefore it is not difficult to see that $\I_{\T}=\I_{t,r}$. 
Formula~\Ref{eq_decomp2_appendix} together with Lemma~\ref{lemma2_appendix} (see below) concludes the proof. 
\end{proof_without_dot}

\begin{lemma}\label{lemma2_appendix} We use notations from Lemma~\ref{lemma1_appendix} and its proof. For any $(\un{j},\un{c})\in \I_{\T}$ we have 
$$\sign(\xi_{\un{j},\un{c}})=(-1)^{\xi},$$ 
where $\xi=t+\sum_{i=1}^{\#\un{c}} j_i(\deg_y{c_i}+\deg_z{c_i}+1)$. Here $c_i$ is considered as a path in $\Q$; therefore $\deg_y{c_i}$ and $\deg_z{c_i}$ are well defined.
\end{lemma}
\begin{proof}
We start the proof with the definition. For a permutation $\pi\in S_{t+2r}$
and $a=a_1\cdots a_s\in \T^{\pi}_{\rm cl}$, where
$a_1,\ldots,a_s$ are arrows or transpose arrows of $\T$, there is a permutation $\tau\in S_{t+2r}$ such that
\begin{enumerate}
\item[$\bullet$] $\tau(i)=i$ for any $i\in
\{1,\ldots,t+2r\}\backslash\{{}'a_1,\ldots,{}'a_s\}$;

\item[$\bullet$] for any
$1\leq i\leq s$ there exists a $b\in\T$ such that $a_i^{\tau}$ and $b$ coincides, i.e., $(a_i^{\tau})''=b''$, $(a_i^{\tau})'=b'$, $''a_i^{\tau}={}''b$, $'a_i^{\tau}={}'b$.
\end{enumerate}
Denote $\sign(a)=\sign(\tau)$. In spite of non-uniqueness of $\tau$, $\sign(a)$ is well defined, it does not depend on $\pi$, and $\sign(a)=\sign(a^T)$. For any $1\leq i\leq \#\un{j}$ we consider $a_i\in\T_{\rm cl}^{\xi_{\un{j},\un{c}}}$ such that $\varphi(a_i)=c_i$. Then we have %
\begin{eq}\label{eq_sgn_xi}
\sign(\xi_{\un{j},\un{c}})=\prod_{i=1}^{\#\un{j}} \sign(a_i)^{j_i}. %
\end{eq} %

Arrows of the tableau $\T$ are $x_i,y_j,z_j$, where $1\leq i\leq t$ and $1\leq j\leq r$. 
Let $\pi\in S_{t+2r}$, $a\in\T_{\rm cl}^{\pi}$, and let $b$, $c$ be some paths in
$\T^{\pi}$.  %
\begin{enumerate}
\item[1.] If $a=x_i$, then $\sign(a)=1$.

\item[2.] If $a=bx_ic$, then
$\sign(a)=-\sign(x_i)\sign(bc)=-\sign(bc)$.

\item[3.] If $a=y_jz_k$, then $\sign(a)=-1$.

\item[4.] If $a=y_jz_k^T$, then $\sign(a)=1$.

\item[5.] If $a=by_jz_kc$, then
$\sign(a)=-\sign(y_jz_k)\sign(bc)=\sign(bc)$.

\item[6.] If $a=by_jz_k^Tc$, then
$\sign(a)=-\sign(y_jz_k^T)\sign(bc)=-\sign(bc)$.

\item[7.] If $a=bz_ky_jc$, then
$\sign(a)=-\sign(y_jz_k)\sign(bc)=\sign(bc)$.

\item[8.] If $a=bz_k^Ty_jc$, then
$\sign(a)=-\sign(y_jz_k^T)\sign(bc)=-\sign(bc)$.
\end{enumerate}

For $a\in \T_{\rm cl}^{\pi}$ we have two possibilities:
\begin{enumerate}
\item[$\bullet$] If $\deg_y a+\deg_z a=0$, then
$\sign(a)=(-1)^{\deg_{x}a+\deg_{x^T}a+1}$ (see items~1 and 2).

\item[$\bullet$] If $\deg_y a+\deg_z a>0$, then by item~2 we have %
$$\sign(a)=(-1)^{\deg_{x}a+\deg_{x^T}a} \sign(a_0),$$ %
where $a_0=a|_{x_i\to 1,\; x_i^T\to 1}$ is the result of
elimination of
letters $x_i,x_i^T (1\leq i\leq t)$ from $a$. By items~$3$--$8$, %
$\sign(a_0)=(-1)^{\deg_{y^T}a_0+\deg_{z^T}a_0+1}= (-1)^{\deg_{z^T}a+\deg_{z^T}a+1}$. 
\end{enumerate}
Thus,  %
$$\sign(a)=(-1)^{\deg_{x}a+\deg_{x^T}a+\deg_{y^T}a+ \deg_{z^T}a+1}.$$
Formula~\Ref{eq_sgn_xi} concludes the proof.
\end{proof}

\begin{lemma}\label{lemma3_appendix}%
Assume that $(\T,(X_1,\ldots,X_t,Y_1,\ldots,Y_r,Z_1\ldots,Z_r))$ is the tableau with substitution from Example~\ref{ex1b_appendix} and $x_1,\ldots,x_t$, $y_1,\ldots,y_r$, $z_1\ldots,z_r$ belong to $\M_{\FF}$. Then 
$$\silin{t,r}(x_1,\ldots,x_t,y_1,\ldots,y_r,z_1\ldots,z_r)=\sum_{\xi\in S_{t+2r}} \sign(\xi) \prod_{a\in \T^{\xi}_{\rm cl}}\tr(\varphi(a))|_{i\to x_i,\;(t+j)\to y_j,\;(t+r+j)\to z_j},$$
where the substitution is applied for all $1\leq i\leq t$ and $1\leq j\leq t$.
In particular, $\silin{t,r}(X_1,\ldots,X_t,Y_1,\ldots,Y_r,Z_1\ldots,Z_r)=
\bpf_{\T}(X_1,\ldots,X_t,Y_1,\ldots,Y_r,Z_1\ldots,Z_r)$ for $n=t+2r$. 
\end{lemma}
\begin{proof} Using the fact that $\T_{\rm cl}^{\xi}\sim \T_{\rm cl}^{\tau}$ for $\xi,\tau\in S_{t+2r}$ if and only if $\xi=\tau$, we prove this lemma similar to Lemma~\ref{lemma1_appendix}.
\end{proof}

\remark\label{remark1_appendix}
Lemma~\ref{lemma2_appendix} can also be generalized for  $\si_{\un{t},\un{r},\un{s}}(x_1,\ldots,x_u,y_1,\ldots,y_v,z_1,\ldots,z_w)$, where $x_1,\ldots,z_w\in\M_{\FF}$.
\bigskip

We assume that $\si'_{t,r}(x,y,z)$ stands for the element of $\N_{\si}$ defined by Zubkov (see p.~292 of~\cite{ZubkovII}), where $t\geq 2r$ and $x,y,z\in\M_{\FF}$.

\begin{lemma}\label{lemmaZ_appendix} For $t,r\geq0$ and $x,y,z\in\M_{\FF}$ we have $\si_{t,r}(x,y,z)=\si'_{t+2r,r}(x,z,y)$.
\end{lemma} 
\begin{proof} 
By Remark~\ref{remark_reduction_to_char_0}, we can assume that $\FF=\QQ$. 
Consider the tableau with substitution $(\T,(X_1,\ldots,X_t,Y_1,\ldots,Y_r,Z_1\ldots,Z_r))$ from Example~\ref{ex1b_appendix}. 
By definition, we have 
$$\si'_{t+2r,r}(x,z,y)=\frac{1}{t!\,(r!)^2}\sum_{\si\in S_{t+2r}} \sign(\si)\, f(\tr\nolimits^{\ast}(\si)),$$
where $\tr^{\ast}(\si)\in\N_{\si}^{\infty}$ is defined on p.~291 of~\cite{ZubkovII} and $f$ stands for the substitution $i\to x_i$, $(t+j)\to y_j$, $(t+r+j)\to z_j$ for all $1\leq i\leq t$ and $1\leq j\leq t$. 
Considering all 18 possibilities from the definition of $\tr^{\ast}(\si)$ we can see that 
$$\tr\nolimits^{\ast}(\si)=\prod_{a\in\T^{\si^{-1}}} \tr(\varphi(a)).$$
Lemmas~\ref{lemma2a_si} and~\ref{lemma3_appendix} conclude the proof.
\end{proof}

\section*{Acknowledgements}
This paper was completed during author's visit to University of Bielefeld, sponsored by DFG LO 1582/1-1. The author is grateful for this support. This research has also been supported by RFFI 10-01-00383a. 



\end{document}